\documentclass[preprint]{elsarticle}
\journal{}

\usepackage{amsmath,amssymb,amsthm}
\usepackage{mathtools}
\usepackage{hyperref}
\usepackage{enumitem}

\newtheorem{theorem}{Theorem}[section]
\newtheorem{lemma}[theorem]{Lemma}
\newtheorem{proposition}[theorem]{Proposition}
\newtheorem{corollary}[theorem]{Corollary}
\newtheorem{conjecture}[theorem]{Conjecture}

\theoremstyle{definition}
\newtheorem{definition}[theorem]{Definition}

\theoremstyle{remark}
\newtheorem{remark}[theorem]{Remark}
\newtheorem{example}[theorem]{Example}

\newcommand{\N}{\mathbb{N}}
\newcommand{\Z}{\mathbb{Z}}
\newcommand{\Q}{\mathbb{Q}}

\newcommand{\Ulam}{\mathrm{U}}     % Ulam predicate U(a,b,m)
\newcommand{\Th}{\mathrm{Th}}      % Theory
\newcommand{\Str}[1]{\mathcal{#1}} % Structures, e.g. \Str{U}

\begin{document}
	
	\begin{frontmatter}
		
		\title{Arithmetical Complexity and Forcing Absoluteness of Rigidity Phenomena for Ulam Sequences}
		
		\author{Frank Gilson\corref{cor1}}
		\cortext[cor1]{Corresponding author.}
		\address{Independent Researcher, Yorba Linda, California}
		
		\begin{abstract}
			We analyse the arithmetical complexity and forcing absoluteness of natural
			statements about Ulam sequences.  For positive integers $a<b$, membership in
			$U(a,b)$ and the increasing enumeration of $U(a,b)$ are uniformly primitive
			recursive.  We encode the finite interval-with-periodic-mask descriptions
			used in rigidity results and formalise the finite-window rigidity theorem for
			the family $U(1,n)$: for every prescribed linear window, one finite collection
			of residue-class-dependent pattern data works for all sufficiently large
			parameters $n$.  This family statement has a $\Pi^0_3$ upper bound in the
			arithmetical hierarchy; for a fixed window its complexity is $\Sigma^0_2$.
			
			For an individual $U(a,b)$, eventual periodicity of the gap sequence (with
			positive period) and finiteness of specified residue classes are $\Sigma^0_2$,
			while rational upper-density, lower-density, and exact-density assertions have
			$\Pi^0_3$ upper bounds.  These are classifications by upper bounds, not
			completeness claims.  Since the resulting sentences are arithmetical, their
			truth is unchanged by set forcing.  This semantic forcing invariance is
			distinguished from proof-theoretic conservativity over stronger
			set theories.  Finally, if the gaps of $U(a,b)$ are eventually periodic with
			positive period, then $U(a,b)$ is Presburger-definable; hence
			$(\N,+,\Ulam_{a,b})$ is decidable, NIP, dp-minimal, and does not interpret
			full arithmetic.
		\end{abstract}
		
		
	\end{frontmatter}
	
	% ------------------------------
	% 1. Introduction
	% ------------------------------
	
	\section{Introduction}
	
	\subsection{Ulam sequences and rigidity}
	
	In his 1964 paper \cite{Ulam1964}, Ulam introduced what is now called the
	\emph{Ulam sequence} $U(1,2)$ as a simple experiment in additive number theory:
	start with $1$ and $2$, and then recursively add the smallest integer that can
	be written as a sum of two distinct earlier terms in a unique way. Thus
	\[
	1,2,3,4,6,8,11,13,16,18,26,28,\dots
	\]
	is obtained by admitting $3=1+2$, then $4=1+3$, while skipping $5$, which can
	be written both as $1+4$ and $2+3$. More generally, for integers $a<b$ one
	defines the Ulam sequence $U(a,b)$ by starting from $\{a,b\}$ and enforcing the
	same ``unique representation as a sum of two distinct earlier terms'' rule.
	
	Despite the simplicity of this definition, Ulam sequences have proved
	remarkably resistant to analysis. Even for the classical sequence $U(1,2)$,
	fundamental basic questions remain open: for example, it is not known whether
	$U(1,2)$ has a natural density, or whether its gap sequence is eventually
	periodic. A series of works by Finch, Schmerl--Spiegel and
	Cassaigne--Finch~\cite{Finch1992ZeroAdditive,SchmerlSpiegel1994,CassaigneFinch1995}
	initiated a systematic study of more general ``$1$-additive'' sequences,
	including $U(a,b)$, and identified many examples in which the gap sequence
	eventually becomes periodic. In particular, they showed that for a number of
	pairs $(a,b)$ the sequence $U(a,b)$ has only finitely many even terms and is
	therefore \emph{regular} in this sense.
	
	A second strand of work has highlighted striking global patterns in Ulam
	sequences. Steinerberger~\cite{Steinerberger2017HiddenSignal}
	observed that for a certain explicit real number $\alpha$, the sequence of
	fractional parts $\{u_n\alpha\}$, where $(u_n)$ enumerates $U(1,2)$, exhibits a
	highly non-uniform limiting distribution concentrated in the middle third of
	$[0,1]$. This ``hidden signal'' is especially surprising in light of Weyl's
	equidistribution theorem, which asserts that for any sequence of distinct
	integers $(a_n)$ and Lebesgue-almost every $\alpha$, the sequence
	$\{a_n\alpha\}$ is equidistributed.
	
	Building on these observations, Hinman, Kuca, Schlesinger and
	Sheydvasser~\cite{HinmanKucaSchlesingerSheydvasser2019JNT,HinmanKucaSchlesingerSheydvasser2019Involve}
	developed a series of \emph{rigidity} results for Ulam sequences.
	Roughly speaking, they showed that as one varies the second term $n$ in $U(1,n)$, 
	large initial segments of the resulting sequences can be described by a finite 
	amount of data that does not depend on $n$ itself, but only on its residue class 
	modulo some fixed $L$. More concretely, for the family $U(1,n)$ they produce 
	decompositions of $U(1,n)\cap[1,c_j n + d_j]$ (where $c_j\to\infty$), valid for 
	all sufficiently large $n$ in a fixed congruence class, as finite unions of 
	intervals with endpoints depending linearly on $n$.
	
	\begin{remark}[The parameter being varied]
		\label{rem:a-greater-1}
		The finite-window theorem of
		Hinman--Kuca--Schlesinger--Sheydvasser fixes the first parameter $a$ and
		varies the second parameter $n$ through residue classes.  We state the
		complexity calculation for $U(1,n)$ to match the best-known presentation;
		the same upper-bound argument applies with any fixed positive $a$.
		Nothing in that theorem asserts that a single fixed $U(a,b)$ is eventually
		equal to a finite union of bounded intervals.
	\end{remark}
	
	These results support the view that Ulam sequences, while seemingly irregular 
	at small scales, exhibit a surprising degree of large-scale order---at least 
	in the $U(1,n)$ family.
	
	\subsection{Main results and significance}
	
	The central contribution of this paper is to determine the \emph{logical 
		complexity} of statements about Ulam sequences and to draw consequences 
	from this classification. Our main theorem, stated informally, is:
	
	\begin{quote}
		\textbf{Main Theorem} (Informal). \emph{Finite-window rigidity for the
			family $U(1,n)$ has a $\Pi^0_3$ arithmetical formulation, while, for a
			fixed sequence $U(a,b)$, positive-period gap regularity is $\Sigma^0_2$
			and rational density assertions have $\Pi^0_3$ upper bounds.  All these
			statements are invariant under set forcing.  If $U(a,b)$ has eventually
			periodic gaps, then the expansion $(\N,+,\Ulam_{a,b})$ is a definitional
			expansion of Presburger arithmetic and is therefore decidable, NIP, and
			dp-minimal.}
	\end{quote}
	
	\subsubsection*{Why this matters: for number theorists}
	
	The rigidity conjectures for Ulam sequences are purely combinatorial questions 
	about sequences of integers. It is natural to ask whether sophisticated 
	set-theoretic machinery might help resolve them---perhaps the answer depends 
	on whether we accept certain axioms, or perhaps forcing techniques could 
	shed light on the problem.
	
	Our results isolate one precise limitation on such methods.  Set forcing does
	not change the truth value of these arithmetical statements, because it does
	not change the standard natural numbers.  This is a semantic invariance
	statement; it does not say that set-theoretic methods cannot prove arithmetic
	statements, nor does it imply conservativity of theories with large-cardinal
	axioms over arithmetic.
	
	More specifically, we show that the key conjectures can be expressed using 
	only two or three alternations of quantifiers over natural numbers. This 
	places them at a very low level of logical complexity---comparable to 
	statements like ``there exist infinitely many twin primes'' or ``Goldbach's 
	conjecture holds.'' The classification gives a precise sense in which these 
	problems, while difficult, are ``elementary'' from a logical standpoint.
	
	\subsubsection*{Why this matters: for model theorists}
	
	Ulam sequences provide a natural family of recursively defined predicates 
	on $\N$ whose model-theoretic behavior is controlled by their combinatorial 
	structure. The expansion $(\N,+,\Ulam_{a,b})$ sits between two well-understood 
	extremes: Presburger arithmetic $(\N,+)$, which is decidable and tame, and 
	full arithmetic $(\N,+,\times)$, which is undecidable and wild.
	
	Under the assumption of combinatorial rigidity, we show that these expansions 
	remain on the tame side: they are definitional expansions of Presburger 
	arithmetic, hence NIP and dp-minimal. Without assuming rigidity, the situation 
	is open and provides concrete test cases for questions about expansions of 
	$(\N,+)$ by sparse additive sets. The paper initiates a systematic study of 
	these expansions and formulates precise conjectures about their tameness.
	
	\subsection{Logical background for non-logicians}
	\label{subsec:logical-background-informal}
	
	For readers unfamiliar with mathematical logic, we provide informal 
	explanations of the key concepts that appear in our results. Precise 
	definitions are given in Section~\ref{sec:preliminaries}.
	
	\subsubsection*{The arithmetical hierarchy}
	
	The \emph{arithmetical hierarchy} classifies statements about natural numbers 
	by the pattern of quantifiers (``for all'' and ``there exists'') needed to 
	express them. The notation $\Sigma^0_n$ and $\Pi^0_n$ indicates $n$ 
	alternating blocks of quantifiers:
	\begin{center}
		\begin{tabular}{ll}
			$\Sigma^0_1$: & $\exists x\, \varphi(x)$ where $\varphi$ is ``directly checkable'' \\
			$\Pi^0_1$: & $\forall x\, \varphi(x)$ where $\varphi$ is ``directly checkable'' \\
			$\Sigma^0_2$: & $\exists x\, \forall y\, \varphi(x,y)$ \\
			$\Pi^0_2$: & $\forall x\, \exists y\, \varphi(x,y)$ \\
			$\Sigma^0_3$: & $\exists x\, \forall y\, \exists z\, \varphi(x,y,z)$ \\
		\end{tabular}
	\end{center}
	Lower levels correspond to ``simpler'' statements. For example, ``$U(1,2)$ 
	has infinitely many even elements'' is $\Pi^0_2$ (for all $N$, there exists 
	an even $m > N$ in $U(1,2)$), while ``the gap sequence of $U(a,b)$ is 
	eventually periodic'' is $\Sigma^0_2$ (there exist $N$ and $p$ such that 
	for all $k \geq N$, the gaps repeat with period $p$).
	
	\subsubsection*{Absoluteness}
	
	A statement is \emph{absolute} between specified models if those models assign
	it the same truth value.  Models of set theory can disagree about the size of
	the continuum or the existence of large cardinals.  Transitive models (and,
	more generally, $\omega$-models) with the same standard natural numbers agree
	on arithmetical truth because they use the same structure $(\N,+,\times)$.
	
	The practical consequence is that set forcing cannot reverse the truth value
	of an arithmetical conjecture about Ulam sequences.  Claims about what stronger
	set theories can prove require separate proof-theoretic arguments and are not
	consequences of arithmetical absoluteness alone.
	
	\subsubsection*{NIP and tameness}
	
	In model theory, structures are classified by how ``wild'' or ``tame'' their 
	definable sets are. A theory is \emph{NIP} (``not the independence 
	property'') if its definable families of sets cannot encode arbitrary 
	combinatorial patterns---roughly, there is a uniform bound on how 
	complicated the interactions between definable sets can be.
	
	Presburger arithmetic $(\N,+)$ is a paradigm of tameness: NIP, decidable, 
	and with a complete understanding of definable sets (they are ``semilinear''). 
	Full arithmetic $(\N,+,\times)$ is wild: it can encode arbitrary computation, 
	and its definable sets are extremely complicated.
	
	When we say that $(\N,+,\Ulam_{a,b})$ is ``tame'' under rigidity assumptions, 
	we mean it inherits the good behavior of Presburger arithmetic rather than 
	the pathologies of full arithmetic.
	
	\subsection{Technical overview}
	
	The starting point of this paper is the observation that all of the objects
	above---the individual Ulam sequences $U(a,b)$, the families $\{U(a,b)\}$, and
	the associated rigidity patterns---are inherently \emph{arithmetical}. For
	each coprime pair $a<b$, membership in $U(a,b)$ is algorithmically decidable:
	given $(a,b,m)$, one can simulate the Ulam construction up to $m$ and check
	whether $m$ is admitted. In particular, each $U(a,b)$ is a recursive subset of
	$\N$, and there is a uniform $\Delta^0_1$ predicate $\Ulam(a,b,m)$ expressing
	$m\in U(a,b)$ in first-order arithmetic. Likewise, the ``finite patterns''
	appearing in rigidity conjectures can be encoded by natural numbers and
	decoded by primitive recursive functions.
	
	From this perspective, Ulam sequences give rise to a large family of first-order
	structures
	\[
	\Str{U}_{a,b} \;=\; (\N,+,0,1,\Ulam_{a,b}),
	\]
	where $\Ulam_{a,b}(m)$ abbreviates $\Ulam(a,b,m)$. Natural questions about
	rigidity, regularity and density then become questions about the first-order
	theory of these expansions and the logical form of particular sentences in the
	language of arithmetic. For example:
	\begin{itemize}[label=\textbullet]
		\item the assertion that the family $U(1,n)$ is uniformly described on every
		prescribed linear window by residue-class-dependent finite patterns becomes
		a sentence with a specific pattern of first-order quantifiers over~$\N$;
		\item the assertion that the gap sequence of $U(a,b)$ is eventually periodic
		can be written as $\exists N\exists p\geq1\forall k\ge N\,(\dots)$, and thus
		occupies a definite place in the arithmetical hierarchy;
		\item statements asserting the existence and value of the asymptotic density
		of $U(a,b)$ can likewise be analysed in terms of their quantifier
		complexity.
	\end{itemize}
	
	Once one knows where these statements sit in the arithmetical hierarchy, a
	number of consequences follow almost for free. In particular, purely
	arithmetical sentences are absolute between transitive models of $\mathrm{ZFC}$
	with the same natural numbers, and their truth value cannot be changed by
	forcing extensions.  Thus the finite-window rigidity, regularity, and density
	statements considered here cannot acquire a different truth value after set
	forcing.  We make no general claim of proof-theoretic conservativity for
	large-cardinal or other extensions of set theory.
	
	Beyond absoluteness, one can also ask about the model-theoretic behaviour of
	the expansions $\Str{U}_{a,b}$. A priori such expansions need not be tame: for
	a general recursive predicate $P\subseteq\N$, the structure $(\N,+,P)$ may
	interpret arbitrary complicated combinatorics, including multiplication, and
	may have the independence property (IP) in the sense of Shelah\cite{ShelahCT} (for these notions, see Section~\ref{sec:model-theory-background}). However, the
	empirical and rigorous rigidity results for Ulam sequences suggest that the
	specific predicates $\Ulam_{a,b}$ could be much better behaved. One of the
	themes of this paper is that strong combinatorial rigidity of $U(a,b)$ should
	correlate with model-theoretic tameness of $\Str{U}_{a,b}$, for instance with
	the absence of IP or with non-interpretability of $(\N,\times)$.
	
	The aim of this work is therefore twofold. First, we give a uniform
	arithmetisation of Ulam sequences and of the ``finite patterns'' used in
	rigidity results, and we identify the precise arithmetical complexity of a
	range of natural statements about $U(a,b)$: finite-window family rigidity, regularity,
	finiteness of even terms, and various density assertions. Second, we draw out
	the resulting absoluteness consequences and outline some initial
	model-theoretic implications for the structures $\Str{U}_{a,b}$. In this way,
	we place the rigidity phenomena for Ulam sequences within the broader landscape
	of definability, absoluteness and classification theory in mathematical logic.
	
	\subsection{Main results}
	
	We now summarise our main contributions more formally. Throughout, we write
	$\Ulam(a,b,m)$ for the uniform $\Delta^0_1$ predicate coding membership in
	$U(a,b)$, as constructed in Theorem~\ref{thm:uniform-ulam}, and
	$\mathrm{InPattern}(c,a,b,m)$ for the uniform primitive recursive relation
	expressing that $m$ belongs to the pattern encoded by $c$ for parameters
	$(a,b)$.
	
	\paragraph{Uniform arithmetisation.}
	Our first result is a uniform arithmetisation of Ulam sequences and of the
	finite ``interval-with-mask'' patterns that arise in rigidity results. In
	Theorem~\ref{thm:uniform-ulam} we construct a single $\Delta^0_1$ predicate
	$\Ulam(a,b,m)$ such that for all coprime $a<b$ and $m\in\N$,
	\[
	\Ulam(a,b,m) \iff m\in U(a,b).
	\]
	In Section~\ref{subsec:coding-finite-data} we define a primitive recursive
	coding of finite patterns and prove (Lemma~\ref{lem:pattern-primitive-recursive})
	that for each code $c$ with $\mathrm{Pattern}(c)$, the associated set
	\[
	\mathrm{PatternSet}_c(a,b) = \{m\in\N : \mathrm{InPattern}(c,a,b,m)\}
	\]
	is a finite union of linearly parameterised intervals with periodic masks, in
	the sense used in rigidity results for $U(a,b)$.
	
	\paragraph{Complexity of finite-window rigidity.}
	For each positive integer $K$, let $\mathrm{WinRig}(K)$ assert that there are
	a modulus, a threshold, and finitely many pattern codes---one for each residue
	class---which correctly describe $U(1,n)$ for every sufficiently large $n$ on
	the window $[1,Kn]$.  Theorem~\ref{thm:window-rigidity-complexity} shows that
	$\mathrm{WinRig}(K)$ is $\Sigma^0_2$ for fixed $K$, and that
	$\forall K\,\mathrm{WinRig}(K)$ is $\Pi^0_3$.  This is the quantifier structure
	of the finite-window rigidity theorem in the literature; it must not be
	confused with eventual equality of a fixed sequence $U(a,b)$ to a finite
	pattern set.
	
	\paragraph{Complexity of regularity and simple combinatorial properties.}
	Let $U(a,b)=\{u_k : k\in\N\}$ be enumerated in increasing order and set
	$g_k = u_{k+1}-u_k$ for the associated gap sequence. In
	Theorem~\ref{thm:regularity-sigma02} we show that, for each fixed coprime
	$(a,b)$, the natural \emph{regularity} statement
	\[
	\exists N\,\exists p\geq1\,\forall k\ge N\,(g_{k+p}=g_k),
	\]
	asserting eventual periodicity with a positive period, is equivalent
	to a $\Sigma^0_2$ sentence of first-order arithmetic. In
	Section~\ref{subsec:finiteness-evens} we show that other basic properties,
	such as ``$U(a,b)$ has only finitely many even elements'', can also be
	expressed by $\Sigma^0_2$ sentences for each fixed $(a,b)$.
	
	\paragraph{Complexity of density statements.}
	In Section~\ref{sec:regularity-density} we study the asymptotic density of
	$U(a,b)$ via the counting function $C(a,b,n)=|U(a,b)\cap[0,n]|$. For each
	fixed coprime $(a,b)$ and rational $q\in\Q$, we show in
	Theorem~\ref{thm:density-complexity} that the upper- and lower-density
	inequalities
	\[
	\overline{d}(U(a,b))\le q, \qquad \underline{d}(U(a,b))\ge q,
	\]
	and the exact-density statement $d(U(a,b))=q$ are all arithmetical and in fact
	equivalent (for fixed parameters) to $\Pi^0_3$ sentences of arithmetic. In
	particular, all of these statements lie low in the arithmetical hierarchy and
	are absolutely determined by the structure $(\N,+,\times)$.
	
	\paragraph{Absoluteness.}
	Since all of the sentences appearing in the rigidity, regularity and density
	results above are arithmetical, standard absoluteness results immediately
	imply that their truth values do not depend on the ambient set-theoretic
	universe. In Section~\ref{sec:absoluteness} we show that if $M$ and $N$ are
	transitive models of $\mathrm{ZFC}$ with the same natural numbers, then for
	each coprime $(a,b)$ and each rational $q$ the following have the same truth
	value in $M$ and~$N$:
	finite-window family rigidity; regularity and finiteness-of-evens;
	and the density statements described above. In particular, these statements
	cannot be made true in one forcing extension of $\mathrm{ZFC}$ and false in
	another while preserving~$\omega$.
	
	\paragraph{Model-theoretic consequences.}
	Finally, we initiate a model-theoretic analysis of the expansions
	\[
	\Str{U}_{a,b} = (\N,+,0,1,\Ulam_{a,b})
	\]
	of Presburger arithmetic by the predicate naming $U(a,b)$. In
	Section~\ref{sec:model-theory} we show that if the gap sequence of $U(a,b)$ is
	eventually periodic with positive period, then the predicate $\Ulam_{a,b}$ is
	definable in the reduct $(\N,+,0,1)$ and $\Str{U}_{a,b}$ is a definitional
	expansion of Presburger arithmetic. As a consequence
	(Proposition~\ref{prop:no-multiplication}), the theory $\Th(\Str{U}_{a,b})$
	does not interpret $(\N,+,\times)$ and remains model-theoretically tame: it is
	NIP and dp-minimal; equivalently, each fixed definable family has finite
	VC-dimension.
	
	The detailed statements and proofs of these results are given in
	Sections~\ref{sec:arithmetisation}--\ref{sec:model-theory} below.
	We emphasise that, beyond the $U(1,n)$ family where rigidity results are
	known, our theorems concern the \emph{logical complexity and absoluteness}
	of the corresponding rigidity, regularity, and density \emph{schemas}
	for general coprime $(a,b)$; we do not claim these combinatorial
	statements hold in full generality.
	
	% Close with an outline:
	
	\subsection{Organisation of the paper}
	
	Section~\ref{subsec:logical-background-informal} together with the technical overview 
	provide informal and formal background for the main results. 
	Section~\ref{sec:preliminaries} gives precise definitions:
	Ulam sequences, the arithmetical hierarchy, absoluteness, and the coding of 
	finite combinatorial patterns. Section~\ref{sec:arithmetisation} constructs 
	the uniform $\Delta^0_1$ predicate $\Ulam(a,b,m)$ and the pattern predicates 
	used in rigidity statements (Theorem~\ref{thm:uniform-ulam}). 
	Section~\ref{sec:rigidity} analyses the complexity of finite-window family
	rigidity, establishing Theorem~\ref{thm:window-rigidity-complexity}.
	Section~\ref{sec:regularity-density} treats regularity and density, proving
	Theorems~\ref{thm:regularity-sigma02} and~\ref{thm:density-complexity}. 
	Section~\ref{sec:absoluteness} draws absoluteness and proof-theoretic 
	consequences. Section~\ref{sec:model-theory} develops the model-theoretic 
	analysis of $(\N,+,\Ulam_{a,b})$, proving Proposition~\ref{prop:no-multiplication} 
	and formulating tameness conjectures. We conclude in 
	Section~\ref{sec:further-directions} with open problems.
	
	\subsection{Notation and conventions}
	
	We collect here the main notational conventions used throughout the paper.
	
	\begin{center}
		\small
		\setlength{\tabcolsep}{3pt}
		\renewcommand{\arraystretch}{1.1}
		\begin{tabular}{p{0.24\textwidth}p{0.71\textwidth}}
			\hline
			Symbol & Meaning \\
			\hline
			$\N,\Z,\Q$
			& The sets of natural, integer, and rational numbers, respectively. \\[3pt]
			$U(a,b)$
			& The Ulam sequence with initial terms $a<b$, seen as a subset of $\N$. \\[3pt]
			$(u_k)_{k\in\N}$
			& The increasing enumeration of $U(a,b)$, so that $U(a,b)=\{u_k:k\in\N\}$. \\[3pt]
			$g_k$
			& The gap sequence $g_k = u_{k+1}-u_k$ associated to $U(a,b)$. \\[3pt]
			$\Ulam(a,b,m)$
			& A fixed $\Delta^0_1$ predicate in arithmetic coding $m\in U(a,b)$
			(Theorem~\ref{thm:uniform-ulam}). \\[3pt]
			$\Ulam_{a,b}$
			& The unary predicate $m\mapsto \Ulam(a,b,m)$ in the structure
			$\Str{U}_{a,b}$. \\[3pt]
			$\mathrm{Pattern}(c)$
			& Primitive recursive predicate asserting that $c\in\N$ encodes
			a well-formed finite pattern
			(Section~\ref{subsec:coding-finite-data}). \\[3pt]
			$\mathrm{InPattern}(c,a,b,m)$
			& Primitive recursive relation expressing that $m$ lies in
			the pattern determined by $c$ at parameters $(a,b)$. \\[3pt]
			$\mathrm{PatternSet}_c(a,b)$
			& The set $\{m\in\N : \mathrm{InPattern}(c,a,b,m)\}$ associated
			to a pattern code $c$. \\[3pt]
			$B_{\max}(c,a,b)$
			& A primitive recursive bound with
			$\mathrm{PatternSet}_c(a,b)\subseteq\{m : m\le B_{\max}(c,a,b)\}$
			(Lemma~\ref{lem:pattern-upper-bound}). \\[3pt]
			$C(a,b,n)$
			& The counting function $C(a,b,n)=|U(a,b)\cap[0,n]|$
			(Lemma~\ref{lem:counting}). \\[3pt]
			$\overline{d}(X),\underline{d}(X)$
			& Upper and lower asymptotic densities of $X\subseteq\N$
			as in Section~\ref{sec:regularity-density}. \\[3pt]
			$d(X)$
			& The (asymptotic) density of $X$ when $\overline{d}(X)=\underline{d}(X)$. \\[3pt]
			$\Str{U}_{a,b}$
			& The expansion $(\N,+,0,1,\Ulam_{a,b})$ of Presburger arithmetic
			by the predicate naming $U(a,b)$ (Section~\ref{sec:model-theory}). \\[3pt]
			$\Th(\mathcal{M})$
			& The complete first-order theory of a structure $\mathcal{M}$. \\
			\hline
		\end{tabular}
	\end{center}
	
	Unless otherwise specified, we work in the standard model $(\N,+,\times,0,1)$ of
	first-order arithmetic. We freely identify subsets of $\N$ with their
	characteristic predicates, and use the usual notation
	$\Sigma^0_n,\Pi^0_n,\Delta^0_n$ for the levels of the arithmetical hierarchy.
	
	\subsection{Motivations}
	Beyond the specific complexity bounds established here, one of our main
	motivations is conceptual. Ulam sequences provide a concrete and well-studied
	test case in additive combinatorics, and the rigidity phenomena observed by
	Hinman--Kuca--Schlesinger--Sheydvasser exhibit a striking blend of structure
	and apparent irregularity. Our results show that a large class of natural
	questions about this behaviour --- including finite-window rigidity, gap
	regularity, and density --- are in fact of very low logical complexity and
	are absolutely determined by the arithmetic of $(\N,+,\times)$. The
	expansions $(\N,+,\Ulam_{a,b})$ form a new family of test structures for the
	model theory of expansions of Presburger arithmetic, linking combinatorial
	rigidity to tameness notions such as NIP and dp-minimality.
	
	% ------------------------------
	% 2. Preliminaries
	% ------------------------------
	
	\section{Preliminaries}\label{sec:preliminaries}
	
	\subsection{Ulam sequences}
	
	We begin by recalling the standard definition of Ulam sequences $U(a,b)$ and
	their higher-dimensional generalisations.
	
	\begin{definition}
		Let $a,b$ be positive integers with $a<b$. The \emph{Ulam sequence} $U(a,b)$ is
		the increasing sequence of positive integers $(u_k)_{k\in\mathbb{N}}$ defined
		recursively as follows:
		\begin{itemize}[label=\textbullet]
			\item $u_1 = a$, $u_2 = b$;
			\item for $k\ge 3$, we set
			\[
			u_k \;=\;
			\min\bigl\{ n>u_{k-1} :
			\exists!\, 1\le i<j<k \,\bigl(n = u_i + u_j\bigr)
			\bigr\}.
			\]
		\end{itemize}
		We write $U(a,b)$ for the underlying set $\{u_k : k\in\mathbb{N}\}$, and we
		write $m\in U(a,b)$ if $m$ is one of these terms.
	\end{definition}
	
	In this paper we will mostly restrict attention to coprime pairs $a<b$; when
	$\gcd(a,b)>1$ the sequence exhibits obvious degeneracies (for instance, it
	remains contained in the arithmetic progression $a+\gcd(a,b)\mathbb{N}$). The
	classical example is $U(1,2)$, introduced by Ulam~\cite{Ulam1964}, whose
	initial terms are
	\[
	U(1,2) = \{1,2,3,4,6,8,11,13,16,18,26,28,\dots\}.
	\]
	More generally, Ulam sequences $U(a,b)$ are special instances of the
	``$1$-additive'' sequences studied by Finch, Schmerl--Spiegel and
	Cassaigne--Finch~\cite{Finch1992ZeroAdditive,SchmerlSpiegel1994,CassaigneFinch1995},
	where each new term is required to have a unique representation as a sum of two
	distinct earlier terms.
	
	Each $U(a,b)$ is infinite.  Indeed, after terms
	$u_1<\cdots<u_k$ have been chosen, the sum $u_{k-1}+u_k$ has exactly one
	representation by two distinct chosen terms: any representation using $u_k$
	forces the other summand to be $u_{k-1}$, while a representation not using
	$u_k$ has smaller sum.  Hence the next term exists and is at most
	$u_{k-1}+u_k$.  This bound also makes the membership relation decidable by a
	bounded simulation of the recursive construction.  We exploit this
	effectivity in Section~\ref{sec:arithmetisation}.
	
	\begin{remark}[Higher-dimensional generalisations]
		Ulam-type constructions can be defined in higher dimensions: 
		Kravitz--Steinerberger~\cite{KravitzSteinerberger2018} introduced 
		\emph{Ulam sets} in $\mathbb{Z}^n$, and 
		Hinman--Kuca--Schlesinger--Sheydvasser~\cite{HinmanKucaSchlesingerSheydvasser2019Involve} 
		established rigidity phenomena for these as well. The arithmetisation 
		methods of this paper extend naturally to such settings, but we focus 
		exclusively on the one-dimensional sequences $U(a,b)$ for concreteness.
	\end{remark}
	
	\subsection{The arithmetical hierarchy and absoluteness}
	
	We briefly recall the arithmetical hierarchy and the basic absoluteness facts
	we will use. Our reference for background material is, for example,
	Shoenfield~\cite{Shoenfield1967} or Kunen~\cite{Kunen1980}.
	
	\medskip
	
	We work in the usual first-order language of arithmetic
	$\mathcal{L}_{\mathrm{PA}} = \{0,1,+,\times,\leq\}$. A formula is
	\emph{bounded} (or $\Delta^0_0$) if every quantifier is of the form
	$\forall x < t$ or $\exists x < t$, where $t$ is a term not containing~$x$.
	The classes $\Sigma^0_n$ and $\Pi^0_n$ of the \emph{arithmetical hierarchy}
	are defined inductively as follows:
	
	\begin{itemize}[label=\textbullet]
		\item $\Sigma^0_0 = \Pi^0_0$ is the class of bounded formulas;
		\item $\Sigma^0_{n+1}$ consists of formulas equivalent (over Peano
		arithmetic) to formulas of the form $\exists x\,\varphi(x,\bar y)$
		with $\varphi\in\Pi^0_n$;
		\item $\Pi^0_{n+1}$ consists of formulas equivalent to formulas of the form
		$\forall x\,\varphi(x,\bar y)$ with $\varphi\in\Sigma^0_n$.
	\end{itemize}
	
	A formula is \emph{arithmetical} if it is $\Sigma^0_n$ or $\Pi^0_n$ for some
	$n$. A set $X\subseteq\N^k$ is said to be $\Sigma^0_n$ (respectively,
	$\Pi^0_n$) if it is definable in $(\N,+,\times)$ by a $\Sigma^0_n$
	(respectively, $\Pi^0_n$) formula; it is $\Delta^0_n$ if it is both
	$\Sigma^0_n$ and $\Pi^0_n$. In particular, \emph{recursive} (computable)
	subsets of $\N^k$ are exactly the $\Delta^0_1$ subsets, and
	\emph{primitive recursive} relations are also $\Delta^0_1$.
	
	Throughout the paper we will freely switch between recursive/primitive
	recursive definitions and low-level arithmetical formulas: whenever we say
	that a relation $R(\bar x)$ is primitive recursive (or recursive), we may
	equivalently view it as being defined by a fixed $\Delta^0_1$ formula
	$R(\bar x)$ in $\mathcal{L}_{\mathrm{PA}}$.
	
	\medskip
	
	The complexity classification in the arithmetical hierarchy interacts well with
	the set-theoretic notion of absoluteness. If $M$ is a transitive model of
	$\mathrm{ZFC}$, we write $\N^M$ for its natural numbers. A sentence
	$\varphi$ in the language of arithmetic is said to be \emph{arithmetical} if
	it is $\Sigma^0_n$ or $\Pi^0_n$ for some $n$ (equivalently, if all of its
	quantifiers range over~$\N$ and it uses only $0,1,+,\times,\le$).
	
	\begin{lemma}[Absoluteness for arithmetical sentences]
		\label{lem:arithmetical-absoluteness}
		Let $M$ and $N$ be transitive models of $\mathrm{ZFC}$ with the same natural
		numbers, i.e.\ $\N^M = \N^N = \N$. Then for every arithmetical sentence
		$\varphi$ in the language of arithmetic, we have
		\[
		M\models\varphi \quad\Longleftrightarrow\quad N\models\varphi.
		\]
	\end{lemma}
	
	\begin{proof}
			The satisfaction relation for arithmetical formulas can be defined by
			primitive recursion on the complexity of formulas, using only
			quantification over~$\N$. In particular, the truth value of an
			arithmetical sentence depends only on the underlying structure
			$(\N,+,\times,\leq)$. Since $M$ and $N$ have the same natural numbers and
			the same interpretations of $+, \times, \leq$, they assign the same truth
			value to every arithmetical sentence.
	\end{proof}
	
	In particular, if $V$ is the ambient universe and $V[G]$ is a forcing
	extension, then $V$ and $V[G]$ have the same natural numbers and the same
	structure $(\N,+,\times,\leq)$, so they agree on all arithmetical sentences.
	Thus the truth value of any arithmetical statement is invariant under set
	forcing. The main complexity results of this paper (for rigidity, regularity
	and density statements about $U(a,b)$) will show that all of these can be
	expressed by arithmetical sentences of low complexity, and hence enjoy this
	strong form of absoluteness.
	
	\subsection{Coding finite data by natural numbers}
	\label{subsec:coding-finite-data}
	
	In this subsection we fix, once and for all, a standard coding of finite
	combinatorial data by natural numbers. This will allow us to treat the
	``finite patterns'' appearing in rigidity statements as single natural numbers,
	and to work uniformly with a primitive recursive decoding procedure.
	
	We assume familiarity with classical G\"odel coding of finite tuples and finite
	subsets of~$\N$. For definiteness, one may proceed as follows.
	
	\begin{itemize}[label=\textbullet]
		\item A finite tuple $(x_0,\dots,x_{k-1})$ of natural numbers is coded by
		\[
		\langle x_0,\dots,x_{k-1}\rangle \;=\;
		p_0^{x_0+1}\cdots p_{k-1}^{x_{k-1}+1},
		\]
		where $(p_i)_{i\in\N}$ is the increasing sequence of prime numbers.
		The corresponding decoding functions $\mathrm{len}$ and
		$\mathrm{entry}(c,i)$, which recover the length and $i$th entry of a
		coded tuple $c$, are primitive recursive.
		\item A finite subset $S\subseteq\{0,\dots,L-1\}$ is coded by the natural
		number
		\[
		\mathrm{code}_L(S) \;=\; \sum_{i\in S} 2^i.
		\]
		Conversely, given $L$ and a code $s<2^L$, membership $i\in S$ is
		equivalent to the $i$th binary digit of $s$ being~$1$, which is again a
		primitive recursive relation.
	\end{itemize}
	
	Signed integer coefficients are coded, for example, by sending $z\in\Z$ to
	$2z$ when $z\geq0$ and to $-2z-1$ when $z<0$.  We will not need to commit to
	a specific coding scheme; any fixed primitive recursive encoding of finite
	tuples and finite subsets will do. What matters is
	that the encoding and decoding operations are primitive recursive, so that all
	of the relations we define below are $\Delta^0_1$.
	
	\medskip
	
	Finite-window rigidity statements for parameterised Ulam families describe
	$U(a,b)$, on a window whose endpoint grows linearly with the varying
	parameters, by a finite union of intervals with linear endpoints and periodic
	masks.  We now formalise the finite data used in such descriptions.
	
	\begin{definition}
		A \emph{pattern component} consists of the following data:
		\begin{itemize}[label=\textbullet]
			\item integers $A_1,A_2,B_1,B_2,p,q,L$ with $L\ge 1$;
			\item a finite subset $S\subseteq\{0,\dots,L-1\}$.
		\end{itemize}
		Given such a component and a pair $(a,b)\in\N^2$, we define the associated
		interval
		\[
		I(a,b) = [A(a,b),B(a,b)]\cap\N,
		\]
		where $A(a,b) = A_1 a + A_2 b + p$ and $B(a,b) = B_1 a + B_2 b + q$, and the
		associated \emph{masked subset}
		\[
		C(a,b) = \bigl\{ m\in I(a,b) : (m-A(a,b)) \bmod L \in S\bigr\}.
		\]
	\end{definition}
	
	A \emph{finite pattern} for a family $U(a,b)$ is then a finite family of
	pattern components, which for each $(a,b)$ gives rise to a finite union of
	masked subsets $C(a,b)$. We encode such patterns by single natural numbers.
	
	\begin{definition}
		A \emph{pattern code} is a natural number $c$ which, under a fixed primitive
		recursive decoding scheme, encodes a finite family of pattern components, i.e.\
		a finite list of tuples of the form
		\[
		(A_1,A_2,B_1,B_2,p,q,L,s),
		\]
		where $L\ge 1$ and $s$ is the code of a subset $S\subseteq\{0,\dots,L-1\}$.
		Given such a code $c$ and parameters $(a,b)\in\N^2$, the associated
		\emph{pattern set} is
		\[
		\mathrm{PatternSet}_c(a,b)
		\;=\; \bigcup_{j<k} C_j(a,b),
		\]
		where $(C_j(a,b))_{j<k}$ ranges over the masked subsets determined by the
		components decoded from~$c$.
	\end{definition}
	
	Formally, we will not work directly with the sets $C_j(a,b)$ and
	$\mathrm{PatternSet}_c(a,b)$, but rather with the corresponding membership
	relations, which we now record.
	
	\begin{lemma}
		\label{lem:pattern-primitive-recursive}
		There exists a primitive recursive predicate $\mathrm{Pattern}(c)$ expressing
		that $c$ encodes a well-formed pattern, and a primitive recursive relation
		$\mathrm{InPattern}(c,a,b,m)$ such that whenever $\mathrm{Pattern}(c)$ holds,
		we have
		\[
		\mathrm{InPattern}(c,a,b,m) \quad\Longleftrightarrow\quad
		m\in \mathrm{PatternSet}_c(a,b).
		\]
		In particular, for each fixed $(a,b)$ the set
		$\mathrm{PatternSet}_c(a,b)\subseteq\N$ is $\Delta^0_1$-definable
		uniformly in~$c$.
	\end{lemma}
	
	\begin{proof}
		By construction, a pattern code $c$ consists of a finite list of components,
		each of which is a finite tuple of integers. Using our fixed coding of finite
		tuples, we may assume that $c$ is itself a code of such a list. The predicate
		$\mathrm{Pattern}(c)$ simply asserts, via bounded primitive recursive checks,
		that $c$ decodes to a list of tuples of the required shape and that each
		component satisfies the basic well-formedness conditions (for instance,
		$L\ge 1$ and the mask code $s$ is less than $2^L$). All of these checks are
		primitive recursive in the entries of the decoded tuples.
		
		Given a well-formed code $c$, to decide whether $m\in\mathrm{PatternSet}_c(a,b)$
		we proceed as follows. We first decode $c$ into its finite list of components.
		For each component $(A_1,A_2,B_1,B_2,p,q,L,s)$ we compute $A(a,b)$ and
		$B(a,b)$ and check whether $m\in[A(a,b),B(a,b)]$. If so, we compute
		$r = (m-A(a,b)) \bmod L$ and test whether the $r$th binary digit of $s$ is~$1$,
		that is, whether $2^r$ occurs in the binary expansion of~$s$. We then take the
		disjunction over all components. Since the number of components and all the
		relevant ranges are bounded by primitive recursive functions of $c,a,b,m$, this
		procedure defines a primitive recursive relation
		$\mathrm{InPattern}(c,a,b,m)$.
		
		By construction, whenever $\mathrm{Pattern}(c)$ holds, the above procedure
		agrees with the intended semantic definition of $\mathrm{PatternSet}_c(a,b)$,
		so the displayed equivalence follows. The final claim about $\Delta^0_1$
		definability is immediate from the observation that primitive recursive
		relations correspond to $\Delta^0_1$ formulas in the language of arithmetic.
	\end{proof}
	
	\begin{remark}[The family $U(1,n)$: mask-free decompositions]
		\label{rem:U1n-no-masks}
		For the specific family $U(1,n)$, the pattern components considered above
		specialise to the ``piecewise linear'' descriptions that have been computed
		experimentally in the literature. In particular, Sheydvasser has implemented
		algorithms which, given the numerical data of $U(1,2),\dots,U(1,14)$, recover
		coefficients $c_j, d_j$ (with $c_j\to\infty$ as $j\to\infty$) such that for 
		all sufficiently large $n$ in a fixed residue class modulo some period $L$, 
		one has a decomposition of the form
		\[
		U(1,n) \cap [1,c_j n + d_j]
		\;=\; \bigcup_{i<k}
		\bigl([a_i n + b_i,\, a'_i n + b'_i]\cap\N\bigr);
		\]
		see \cite{Sheydvasser2021LinearPolynomials,SheydvasserUlamCode}. 
		Crucially, for the family $U(1,n)$ there are \emph{no masks}: every interval 
		$[a_i n + b_i,\, a'_i n + b'_i]$ contributes all of its integer points to 
		$U(1,n)$, without periodic deletions. In our notation, this corresponds to 
		pattern components with $L=1$ and $S=\{0\}$ (the trivial mask).
		
		This mask-free property is specific to the $U(1,n)$ family. For general 
		$U(a,b)$ with $a>1$, rigidity phenomena are less well understood, and 
		non-trivial masks may in principle appear; see Remark~\ref{rem:a-greater-1} 
		below. The general pattern-coding framework of 
		Section~\ref{subsec:coding-finite-data} accommodates both cases.
	\end{remark}
	
	In the sequel, we identify a pattern code $c$ with the associated family of
	sets $(\mathrm{PatternSet}_c(a,b))_{(a,b)\in\N^2}$.  The comparison with an
	Ulam sequence will always be restricted to an explicitly bounded window.
	This restriction is essential: for fixed $(c,a,b)$ the pattern set is finite,
	whereas $U(a,b)$ is infinite.
	
	% ------------------------------
	% 3. Uniform coding of U(a,b)
	% ------------------------------
	
	\section{Uniform coding of Ulam sequences in arithmetic}
	\label{sec:arithmetisation}
	
	\subsection{A $\Delta^{0}_{1}$ definition of $\Ulam(a,b,m)$}
	
	We first show that membership in $U(a,b)$ is uniformly decidable by a Turing
	machine which takes $(a,b,m)$ as input. We then invoke the standard fact that
	recursive relations on~$\N$ are $\Delta^0_1$-definable in $(\N,+,\times)$ to
	obtain a uniform $\Delta^0_1$ predicate.
	
	\begin{lemma}[Algorithmic decidability of membership]
		\label{lem:ulam-recursive}
		There is a Turing machine $M$ which, given as input a triple $(a,b,m)$ of
		positive integers with $a<b$, halts with output $1$ if $m\in U(a,b)$ and with
		output $0$ otherwise. In particular, the ternary relation
		\[
		R(a,b,m)\;:\!\iff\; m\in U(a,b)
		\]
		is recursive (computable), uniformly in $(a,b)$.
	\end{lemma}
	
	\begin{proof}
		Fix $a<b$ and $m\in\N$. We describe an algorithm which decides whether
		$m\in U(a,b)$. The algorithm will construct the initial segment of the Ulam
		sequence $U(a,b)$ up to the point where either $m$ has appeared or it can no
		longer appear.
		
		\smallskip
		
		\noindent\emph{Initialisation.}
		Set $S_2 = \{a,b\}$ and $u_1=a$, $u_2=b$. If $m\in\{a,b\}$, halt with output
		$1$.
		
		\smallskip
		
		\noindent\emph{Inductive step.}
		Suppose we have constructed $S_k = \{u_1,\dots,u_k\}$, an initial segment of
		$U(a,b)$, where $u_1<u_2<\dots<u_k$ and $k\ge 2$. We now search for the next
		Ulam term $u_{k+1}$ as follows.
		
		For each integer $n$ strictly larger than $u_k$, in increasing order, we test
		whether $n$ has a \emph{unique} representation as a sum of two distinct
		elements of $S_k$. That is, we count the number of unordered pairs
		$\{i,j\}\subseteq\{1,\dots,k\}$ with $i\neq j$ such that $u_i+u_j=n$. This
		involves only finitely many additions and comparisons, hence is a finite
		computation. We proceed with $n=u_k+1,u_k+2,\dots$, stopping at the first $n$
		for which there is exactly one such representation.
		
		There are two possibilities:
		
		\begin{itemize}[label=\textbullet]
			\item[(1)] If the first such $n$ satisfies $n\le m$, we set $u_{k+1}=n$ and
			$S_{k+1}=S_k\cup\{n\}$, and repeat the inductive step with $k:=k+1$. If
			$n=m$, then at this stage we have found $m$ as an element of $U(a,b)$ and we
			may halt with output~$1$.
			\item[(2)] If the first such $n$ satisfies $n>m$, then no further Ulam terms
			are $\le m$. Indeed, by definition $u_{k+1}$ is the smallest admissible
			integer larger than $u_k$, so any subsequent term $u_{k'}$ (for $k'>k$) must
			satisfy $u_{k'}\ge u_{k+1}>m$. In particular, no new elements of
			$U(a,b)\cap[1,m]$ will ever appear beyond this stage. If we reach such an
			$n>m$ without having seen $m$ among $\{u_1,\dots,u_k\}$, then
			$m\notin U(a,b)$ and we may halt with output~$0$.
		\end{itemize}
		
		\noindent The search is bounded: $u_{k-1}+u_k$ has a unique representation
		as a sum of two distinct elements of $S_k$, so
		$u_{k+1}\leq u_{k-1}+u_k$.  Thus there is always a next term and each search
		terminates.
		Therefore, for each fixed $(a,b,m)$, the above procedure halts after finitely
		many steps with the correct answer. Moreover, the algorithm is uniform: the
		code of the Turing machine does not depend on $(a,b)$ or~$m$ and operates the
		same way for all inputs.
		
		All loops above have explicit primitive-recursive bounds (there are at most
		$m$ admitted positive integers not exceeding $m$, and each candidate search
		is bounded by the sum of the two largest current terms).  Thus the relation
		is in fact primitive recursive.  In particular, the characteristic function of the relation
		$R(a,b,m)\iff m\in U(a,b)$ is computable. Hence $R$ is a recursive ternary
		relation on~$\N$, uniformly in the parameters $(a,b)$.
	\end{proof}
	
	We now translate this into a $\Delta^0_1$ definition in the language of
	arithmetic. The argument is standard: recursive relations are exactly the
	relations definable by $\Delta^0_1$ formulas in $(\N,+,\times)$, via Kleene's
	$T$-predicate and the representability of recursive functions in Peano
	arithmetic (see, for instance, \cite[Chapter~VI]{Shoenfield1967}).
	
	\begin{theorem}[Uniform arithmetisation of Ulam sequences]
		\label{thm:uniform-ulam}
		There exists a $\Delta^{0}_{1}$ formula $\Ulam(a,b,m)$ in the language of
		arithmetic such that for all coprime $a<b$ and all $m\in\N$, we have
		\[
		\Ulam(a,b,m) \iff m\in U(a,b).
		\]
		Moreover, the definition is uniform in $(a,b)$: the same formula
		$\Ulam(a,b,m)$ works for all input parameters $(a,b)$.
	\end{theorem}
	
	\begin{proof}
		By Lemma~\ref{lem:ulam-recursive}, the relation
		\[
		R(a,b,m)\;:\!\iff\; m\in U(a,b)
		\]
		is recursive (computable) as a subset of $\N^3$, uniformly in $(a,b)$. It is a
		standard theorem of recursion theory that every recursive relation on~$\N$ is
		definable in $(\N,+,\times)$ by a $\Delta^0_1$ formula: concretely, one can
		fix a universal Turing machine, code computation histories by natural numbers,
		and define a $\Sigma^0_1$ formula expressing that a computation halts with a
		given output in a given number of steps; the complement of a recursive
		relation is also recursively enumerable, so the same relation is $\Pi^0_1$ as
		well, hence $\Delta^0_1$ (see, e.g., \cite[Theorem~VI.1.1]{Shoenfield1967}).
		
		Applying this to the recursive relation $R(a,b,m)$, we obtain a
		$\Delta^0_1$ formula $\Ulam(a,b,m)$ such that
		\[
		(\N,+,\times)\models \Ulam(a,b,m) \quad\Longleftrightarrow\quad R(a,b,m),
		\]
		i.e.\ $\Ulam(a,b,m)$ holds exactly when $m\in U(a,b)$. The construction of
		$\Ulam$ depends only on the fixed universal Turing machine and the coding of
		computation histories, not on the particular values of $(a,b)$, so the
		definition is uniform in $(a,b)$ as claimed.
	\end{proof}
	
	\subsection{Coding patterns for rigidity}
	
	We now package the pattern machinery from Section~\ref{subsec:coding-finite-data}
	into a form tailored to rigidity statements. Recall that
	Lemma~\ref{lem:pattern-primitive-recursive} provides a primitive recursive
	predicate $\mathrm{Pattern}(c)$ and a primitive recursive relation
	$\mathrm{InPattern}(c,a,b,m)$ indicating membership of $m$ in the pattern
	determined by $c$ at parameters $(a,b)$.
	
	\begin{definition}
		Given a pattern code $c$ satisfying $\mathrm{Pattern}(c)$, we define the
		associated \emph{pattern set} at parameters $(a,b)\in\N^2$ by
		\[
		\mathrm{PatternSet}_c(a,b) \;=\;
		\{m\in\N : \mathrm{InPattern}(c,a,b,m)\}.
		\]
	\end{definition}
	
	Thus for each fixed $c$ and $(a,b)$, the set $\mathrm{PatternSet}_c(a,b)$ is
	the finite union of masked intervals $C_j(a,b)$ determined by the components
	encoded in~$c$, as described in Section~\ref{subsec:coding-finite-data}. We
	record two elementary facts that will be used repeatedly.
	
	\begin{definition}[Masked interval]
		\label{def:masked-interval}
		Let $A,B,L\in\Z$ with $A\le B$ and $L\ge 1$, and let
		$S\subseteq\{0,\dots,L-1\}$ be finite. The \emph{masked interval} determined by
		$(A,B,L,S)$ is the set
		\[
		I(A,B,L,S) = \{ m\in\N : A\le m\le B \text{ and } (m-A)\bmod L \in S \}.
		\]
		Thus $I(A,B,L,\{0,\dots,L-1\})$ is the full interval $[A,B]\cap\N$, while
		proper subsets $S$ encode periodic deletion of residue classes modulo $L$.
	\end{definition}
	
	\begin{lemma}
		\label{lem:pattern-delta01}
		For each pattern code $c$, the relation
		\[
		P_c(a,b,m)\;:\!\iff\; \mathrm{InPattern}(c,a,b,m)
		\]
		is $\Delta^0_1$ in the structure $(\N,+,\times)$; in particular, for each fixed
		$c$ and $(a,b)$ the set $\mathrm{PatternSet}_c(a,b)\subseteq\N$ is recursive.
	\end{lemma}
	
	\begin{proof}
		By Lemma~\ref{lem:pattern-primitive-recursive}, the relation
		$\mathrm{InPattern}(c,a,b,m)$ is primitive recursive as a subset of
		$\N^4$, uniformly in $c$. Every primitive recursive relation is recursive,
		hence definable by a $\Delta^0_1$ formula in the language of arithmetic (see,
		for instance, \cite[Chapter~VI]{Shoenfield1967}). This yields the claim.
	\end{proof}
	
	Next we isolate an explicit upper bound on the support of a pattern set.
	
	\begin{lemma}
		\label{lem:pattern-upper-bound}
		There exists a primitive recursive function $B_{\max}(c,a,b)$ such that for
		every pattern code $c$ with $\mathrm{Pattern}(c)$ and every $(a,b)\in\N^2$ we
		have
		\[
		\mathrm{PatternSet}_c(a,b) \subseteq \{ m\in\N : m \le B_{\max}(c,a,b)\}.
		\]
		In particular, for fixed $(c,a,b)$ the set $\mathrm{PatternSet}_c(a,b)$ is
		finite, with size bounded by a primitive recursive function of $(c,a,b)$.
	\end{lemma}
	
	\begin{proof}
		Fix a pattern code $c$ with $\mathrm{Pattern}(c)$. By definition, $c$ encodes a
		finite list of components
		\[
		(A_1^{(j)},A_2^{(j)},B_1^{(j)},B_2^{(j)},p^{(j)},q^{(j)},L^{(j)},s^{(j)})
		\quad (j<k),
		\]
		each of which contributes a masked subset $C_j(a,b)\subseteq\N$ contained in an
		interval
		\[
		I_j(a,b) = [A^{(j)}(a,b),B^{(j)}(a,b)]\cap\N,
		\]
		where $A^{(j)}(a,b) = A_1^{(j)} a + A_2^{(j)} b + p^{(j)}$ and
		$B^{(j)}(a,b) = B_1^{(j)} a + B_2^{(j)} b + q^{(j)}$. The pattern set is
		\[
		\mathrm{PatternSet}_c(a,b)
		= \bigcup_{j<k} C_j(a,b)
		\subseteq \bigcup_{j<k} I_j(a,b).
		\]
		
		Define
		\[
		B_{\max}(c,a,b)
		= \max_{j<k} B^{(j)}(a,b).
		\]
		Since $k$ and the tuples $(A_1^{(j)},\dots,q^{(j)})$ are all decoded from $c$
		by primitive recursive functions, the function $B_{\max}$ is primitive
		recursive in $(c,a,b)$. By construction, for every $j<k$ and every
		$m\in I_j(a,b)$, we have $m\le B^{(j)}(a,b)\le B_{\max}(c,a,b)$. Hence every
		$m\in\mathrm{PatternSet}_c(a,b)$ satisfies $m\le B_{\max}(c,a,b)$, as required.
		
		The final statement about finiteness follows immediately, since a subset of
		$\{0,1,\dots,B_{\max}(c,a,b)\}$ is finite and its size is bounded by
		$B_{\max}(c,a,b)+1$, which is primitive recursive in $(c,a,b)$.
	\end{proof}
	
	In Section~\ref{sec:rigidity} we compare $U(1,n)$ with such pattern sets only
	on bounded windows $[1,Kn]$.  Lemma~\ref{lem:pattern-delta01} makes each
	membership test $\Delta^0_1$, and the finite window makes the comparison a
	bounded arithmetical condition once $K$ and $n$ are fixed.
	
	\begin{example}[Encoding a single linear block]
		\label{ex:pattern-coding}
		As a concrete illustration, consider the family $U(1,n)$ for $n\ge 4$ studied in
		\cite{HinmanKucaSchlesingerSheydvasser2019JNT,Sheydvasser2021LinearPolynomials}.
		Table~1 of~\cite{Sheydvasser2021LinearPolynomials} (reproduced there as
		``The first 30 intervals of $U(1,X)$'') shows, among others, a block of the form
		\[
		[4X+2,\,5X-1] \subseteq U(1,X),
		\]
		which upon evaluation at $X=n$ gives the interval
		$[4n+2,\,5n-1]\subseteq U(1,n)$ for all sufficiently large $n$.
		
		In our notation, this single block can be realised by a single pattern component
		with coefficients
		\[
		A_1 = 0,\quad A_2 = 4,\quad
		B_1 = 0,\quad B_2 = 5,\quad
		p = 2,\quad q = -1,\quad
		L = 1,\quad S = \{0\}.
		\]
		For parameters $(a,b)=(1,n)$ this gives
		\[
		A(1,n) = A_1\cdot 1 + A_2\cdot n + p = 4n+2,\qquad
		B(1,n) = B_1\cdot 1 + B_2\cdot n + q = 5n-1,
		\]
		so the associated contribution to $\mathrm{PatternSet}_c(1,n)$ is exactly
		\[
		\{m\in\N : 4n+2 \le m \le 5n-1\}.
		\]
		since the mask $S=\{0\}$ and period $L=1$ simply mean that we take all points in
		the interval. More complicated components, with $L>1$ and $S\subsetneq\{0,\dots,L-1\}$,
		encode periodic deletion of residue classes inside such linearly parameterised
		intervals; see Section~\ref{subsec:coding-finite-data} for the general scheme.
		
		For instance, taking $L=3$ and $S=\{0,2\}$ would keep only those $m$ in the
		interval with $(m-(4n+2))\bmod 3\in\{0,2\}$, that is, deleting every third element.
	\end{example}
	
	% ------------------------------
	% 4. Complexity of rigidity statements
	% ------------------------------
	
	\section{Complexity of finite-window rigidity}
	\label{sec:rigidity}

	The rigidity theorem of
	Hinman--Kuca--Schlesinger--Sheydvasser concerns a parameterised family, not
	the tail of one fixed Ulam sequence.  For each positive integer $K$, it gives
	uniform finite descriptions of $U(a,n)$ on $[1,Kn]$ for all sufficiently
	large $n$, with the description allowed to depend on the residue class of
	$n$ modulo a fixed modulus.  This distinction is essential.  Every pattern
	set from Section~\ref{subsec:coding-finite-data} is finite after the
	parameters are fixed, while every Ulam sequence is infinite; consequently,
	eventual equality of a fixed $U(a,b)$ with one such finite pattern set would
	be impossible.

	We code a modulus $L\geq1$ together with a list
	$\vec c=(c_0,\ldots,c_{L-1})$ of pattern codes by one natural number $d$.
	The predicate $\mathrm{FamilyCode}(d,L)$ and the decoding function
	$\mathrm{CodeAt}(d,r)$ are primitive recursive.  Define the bounded comparison
	predicate
	\[
	\begin{aligned}
	\mathrm{Matches}(d,L,n,K)\quad\Longleftrightarrow\quad
	&\mathrm{FamilyCode}(d,L)\\
	&{}\wedge\ \forall m\leq Kn\,
	\Bigl(1\leq m\ \longrightarrow\\
	&\qquad \bigl(\Ulam(1,n,m)\leftrightarrow\\
	&\qquad\quad
	 \mathrm{InPattern}(\mathrm{CodeAt}(d,n\bmod L),1,n,m)\bigr)\Bigr).
	\end{aligned}
	\]
	Because the quantifier over $m$ is bounded and all decoding and membership predicates
	are $\Delta^0_1$, $\mathrm{Matches}$ is $\Delta^0_1$ uniformly in its
	arguments.

	\begin{definition}[Finite-window rigidity]
		For $K\geq1$, let $\mathrm{WinRig}(K)$ be the statement
		\[
		\exists L\geq1\,\exists N_0\,\exists d\,
		\forall n\geq N_0\;\mathrm{Matches}(d,L,n,K).
		\]
		The full finite-window assertion is
		\[
		\mathrm{FWRig}\quad:\!\Longleftrightarrow\quad
		\forall K\geq1\;\mathrm{WinRig}(K).
		\]
	\end{definition}

	The use of integral $K$ entails no loss for windows $[1,Cn]$ with fixed
	positive rational or real $C$: enlarging $C$ to an integer only enlarges the
	window.  Conversely, the integer instances are among the real instances.

	\begin{theorem}[Complexity of finite-window rigidity]
		\label{thm:window-rigidity-complexity}
		For each fixed positive integer $K$, $\mathrm{WinRig}(K)$ is equivalent
		to a $\Sigma^0_2$ sentence.  The full assertion $\mathrm{FWRig}$ is
		equivalent to a $\Pi^0_3$ sentence.
	\end{theorem}

	\begin{proof}
		For fixed $K$, absorb $L,N_0,d$ into one existentially quantified code.
		Since $\mathrm{Matches}$ is $\Delta^0_1$, the resulting prefix is
		$\exists w\,\forall n$ followed by a $\Delta^0_1$ matrix.  This is a
		$\Sigma^0_2$ upper bound.  Quantifying universally over $K$ yields
		\[
		\forall K\,\exists w\,\forall n\;\Theta(K,w,n),
		\]
		with $\Theta$ $\Delta^0_1$, hence a $\Pi^0_3$ upper bound.
	\end{proof}

	\begin{remark}[Relation to the published rigidity theorem]
		Theorem~2.1 of
		\cite{HinmanKucaSchlesingerSheydvasser2019JNT} proves the corresponding
		finite-window decomposition for every fixed positive first parameter
		$a$.  The theorem above classifies the arithmetical form of that assertion;
		it does not reprove the combinatorial decomposition.  A stronger conjecture
		postulating one infinite sequence of interval coefficients requires a
		separate coding choice.  Unless an effective code for that infinite object
		is included, it should not be identified with the first-order sentence
		$\mathrm{FWRig}$.
	\end{remark}

	\section{Complexity of regularity and density statements}
	\label{sec:regularity-density}
	
	\subsection{Eventual periodicity of gaps}
	
	We begin by formalising the notion of regularity via the gap sequence.
	
	\begin{definition}
		Let $(u_k)_{k\in\N}$ be the increasing enumeration of $U(a,b)$. The
		\emph{gap sequence} associated to $U(a,b)$ is the sequence
		\[
		g_k = u_{k+1} - u_k \qquad (k\in\N).
		\]
		We say that $U(a,b)$ is \emph{regular} if there exist $N,p\in\N$ with
		$p\geq1$ such that
		for all $k\ge N$ we have $g_{k+p}=g_k$, i.e.
		\[
		\exists N\,\exists p\geq1\,\forall k\ge N\,
		\bigl(u_{k+1}-u_k = u_{k+p+1}-u_{k+p}\bigr).
		\]
	\end{definition}
	
	To analyse the logical complexity of this statement, we need to know that the
	enumeration $k\mapsto u_k$ is itself recursive (and hence $\Delta^0_1$-definable)
	uniformly in $(a,b)$. This is immediate from the algorithmic construction of
	$U(a,b)$.
	
	\begin{lemma}[Recursive enumeration of $U(a,b)$]
		\label{lem:ulam-enumeration}
		There is a Turing machine $M'$ which, given as input a triple $(a,b,k)$ of
		positive integers with $a<b$, halts with output $u_k$, where $(u_k)$ is the
		increasing enumeration of $U(a,b)$. In particular, the function
		\[
		E(a,b,k) = u_k
		\]
		is a total recursive function of $(a,b,k)$. Consequently, there is a
		$\Delta^0_1$ formula $\mathrm{UTerm}(a,b,k,m)$ in the language of arithmetic
		such that for all $(a,b,k,m)\in\N^4$,
		\[
		\mathrm{UTerm}(a,b,k,m) \quad\Longleftrightarrow\quad m = u_k.
		\]
	\end{lemma}
	
	\begin{proof}
		Fix $a<b$ and $k\in\N$. We describe an algorithm which computes $u_k$. The
		construction is essentially the same as in the proof of
		Lemma~\ref{lem:ulam-recursive} but now we stop when we reach the $k$th term
		rather than a fixed bound $m$.
		
		\smallskip
		
		\noindent\emph{Initialisation.}
		Set $S_2=\{a,b\}$ and $u_1=a$, $u_2=b$. If $k=1$ (respectively, $k=2$), halt
		with output $a$ (respectively, $b$).
		
		\smallskip
		
		\noindent\emph{Inductive step.}
		Suppose $k\ge 3$ and we have constructed $S_j=\{u_1,\dots,u_j\}$ for some
		$j\ge 2$, where $u_1<\dots<u_j$ are the first $j$ terms of $U(a,b)$. To
		construct $u_{j+1}$, we search over integers $n>u_j$ in increasing order and
		count the number of unordered pairs $\{i,\ell\}\subseteq\{1,\dots,j\}$ with
		$i\neq \ell$ and $u_i+u_\ell=n$. We stop at the first $n$ for which there is
		exactly one such representation, and declare $u_{j+1}=n$. We then set
		$S_{j+1}=S_j\cup\{n\}$ and continue.
		
		This procedure constructs the Ulam sequence $U(a,b)$ term by term. Since
		$U(a,b)$ is infinite and strictly increasing, the process continues
		indefinitely, and for each input $k$ the algorithm reaches the stage $j=k$
		after finitely many steps, at which point $u_k$ has been computed. Thus the
		function $E(a,b,k) = u_k$ is computable (recursive) and uniform in the
		parameters $(a,b,k)$.
		
		As in the proof of Theorem~\ref{thm:uniform-ulam}, every total recursive
		function on~$\N^3$ is representable in $(\N,+,\times)$ by a $\Delta^0_1$
		formula: there is a $\Delta^0_1$ formula $\mathrm{UTerm}(a,b,k,m)$ such that
		$\mathrm{UTerm}(a,b,k,m)$ holds if and only if $m$ is the output of the fixed
		Turing machine $M'$ on input $(a,b,k)$; see, for instance,
		\cite[Theorem~VI.1.1]{Shoenfield1967}. By construction, this output is $u_k$,
		so $\mathrm{UTerm}(a,b,k,m)\leftrightarrow m=u_k$ as required.
	\end{proof}
	
	We can now express regularity entirely in terms of the $\Delta^0_1$ predicate
	$\mathrm{UTerm}$.
	
	\begin{theorem}
		\label{thm:regularity-sigma02}
		For each fixed coprime $(a,b)$, the statement
		\[
		\exists N\,\exists p\geq1\,\forall k\ge N\,\bigl(g_{k+p}=g_k\bigr)
		\]
		expressing regularity of $U(a,b)$ is equivalent to a $\Sigma^0_2$ sentence of
		first-order arithmetic.
	\end{theorem}
	
	\begin{proof}
		By Lemma~\ref{lem:ulam-enumeration}, the function $E(a,b,k)=u_k$ is
		total recursive; using the explicit bound in the construction, it may be
		taken primitive recursive.  Hence
		\[
		\begin{aligned}
		\mathrm{GapEq}(a,b,k,p)\quad\Longleftrightarrow\quad
		&E(a,b,k+1)+E(a,b,k+p)\\
		&=E(a,b,k)+E(a,b,k+p+1).
		\end{aligned}
		\]
		is $\Delta^0_1$ and expresses $g_{k+p}=g_k$.  Regularity is therefore
		equivalent to
		\[
		\exists N\,\exists p\,\forall k\,
		\bigl(p\geq1\ \wedge\ (k<N\ \vee\
		\mathrm{GapEq}(a,b,k,p))\bigr).
		\]
		The matrix is $\Delta^0_1$, so this is a $\Sigma^0_2$ sentence.
	\end{proof}

	\subsection{Finiteness of even terms and related properties}\label{subsec:finiteness-evens}
	
	We next consider simple combinatorial properties of $U(a,b)$ such as having
	only finitely many even elements. These also sit low in the arithmetical
	hierarchy once expressed in terms of the predicate $\Ulam(a,b,m)$.
	
	\begin{proposition}\label{prop:finiteness-evens}
		For fixed coprime $(a,b)$, the statement \emph{``$U(a,b)$ has only finitely
			many even elements''} is equivalent to a $\Sigma^0_2$ sentence of arithmetic.
	\end{proposition}
	
	\begin{proof}
		Fix coprime $a<b$. The set of even elements of $U(a,b)$ is
		\[
		E(a,b) = \{ m\in\N : \text{$m$ is even and } \Ulam(a,b,m)\}.
		\]
		Saying that $E(a,b)$ is finite is equivalent to saying that there exists an
		$N$ such that no even $m\ge N$ belongs to $U(a,b)$, i.e.
		\[
		\exists N\,\forall m\,
		\Bigl(m\ge N \wedge 2\mid m \ \rightarrow\ \neg \Ulam(a,b,m)\Bigr).
		\]
		Call this sentence $(\ast)$.
		
		First, if $(\ast)$ holds, then there is some $N$ such that every even element
		of $U(a,b)$ is less than $N$, so $E(a,b)\subseteq\{0,2,4,\dots,N-2\}$ and is
		in particular finite. Conversely, if $E(a,b)$ is finite, there is some $N$
		greater than all elements of $E(a,b)$. Then by definition no even $m\ge N$
		lies in $U(a,b)$, so $(\ast)$ holds. Thus $(\ast)$ is equivalent (over Peano
		arithmetic) to the informal statement that $U(a,b)$ has only finitely many
		even elements.
		
		We now analyse the logical form of $(\ast)$. The outer quantifiers are
		$\exists N\,\forall m$. The matrix
		\[
		m\ge N \wedge 2\mid m \ \rightarrow\ \neg \Ulam(a,b,m)
		\]
		is built from bounded inequalities, the divisibility relation $2\mid m$ and
		the predicate $\Ulam(a,b,m)$. The inequality $m\ge N$ and the divisibility
		condition $2\mid m$ are $\Delta^0_0$ formulas, and $\Ulam(a,b,m)$ is
		$\Delta^0_1$ by Theorem~\ref{thm:uniform-ulam}. Hence the whole matrix is
		$\Delta^0_1$. Therefore $(\ast)$ has the form
		\[
		\exists N\,\forall m\,\Phi(N,m),
		\]
		with $\Phi$ a $\Delta^0_1$ formula. By the definition of the arithmetical
		hierarchy, this is a $\Sigma^0_2$ sentence. This completes the proof.
	\end{proof}
	
	\begin{remark}
		The negation of the finiteness statement, namely \emph{``$U(a,b)$ has
			infinitely many even elements''}, can be written as
		\[
		\forall N\,\exists m\ge N\,\bigl(2\mid m \wedge \Ulam(a,b,m)\bigr),
		\]
		which has the quantifier prefix $\forall\exists$ with a $\Delta^0_1$ matrix
		and hence is $\Pi^0_2$. Similar reasoning applies to other simple properties
		such as \emph{``$U(a,b)$ contains only finitely many terms in a given residue
			class modulo $\ell$''} or \emph{``$U(a,b)$ contains infinitely many terms in
			a given residue class modulo $\ell$''}, which are again $\Sigma^0_2$ and
		$\Pi^0_2$ respectively for each fixed $(a,b)$ and~$\ell$.
	\end{remark}
	
	\subsection{Density}
	
	We recall the usual notions of upper and lower asymptotic density.
	
	\begin{definition}
		For a set $X\subseteq\N$, the upper and lower asymptotic densities are
		\[
		\overline{d}(X) = \limsup_{n\to\infty}\frac{|X\cap[0,n]|}{n+1}, \qquad
		\underline{d}(X) = \liminf_{n\to\infty}\frac{|X\cap[0,n]|}{n+1}.
		\]
		If these coincide, we write $d(X)$ for the common value.
	\end{definition}
	
	For $U(a,b)$ we will work with the associated counting function
	\[
	C_{a,b}(n) \;=\; |U(a,b)\cap[0,n]|.
	\]
	Since membership in $U(a,b)$ is recursive, so is $C_{a,b}(n)$.
	
	\begin{lemma}[Counting function]
		\label{lem:counting}
		There is a total recursive function $C(a,b,n)$ such that for all
		$(a,b,n)\in\N^3$,
		\[
		C(a,b,n) = |U(a,b)\cap[0,n]|.
		\]
		Consequently, there is a $\Delta^0_1$ formula $\mathrm{Count}(a,b,n,t)$ in the
		language of arithmetic such that
		\[
		\mathrm{Count}(a,b,n,t) \quad\Longleftrightarrow\quad
		t = C(a,b,n).
		\]
	\end{lemma}
	
	\begin{proof}
		Fix $a<b$ and $n\in\N$. To compute $|U(a,b)\cap[0,n]|$, we may proceed as
		follows. Using the algorithm from Lemma~\ref{lem:ulam-recursive}, we can
		effectively decide for each $m\le n$ whether $m\in U(a,b)$, by constructing
		$U(a,b)$ up to the stage at which either we see $m$ or the next Ulam term
		exceeds $m$. Running this decision procedure for all $m=0,1,\dots,n$ and
		incrementing a counter whenever $m\in U(a,b)$ yields the desired cardinality.
		
		This gives a Turing machine which, on input $(a,b,n)$, halts with output
		$|U(a,b)\cap[0,n]|$. Hence the function $C(a,b,n)$ is recursive. As in the
		proof of Theorem~\ref{thm:uniform-ulam}, every total recursive function on
		$\N^3$ is representable in $(\N,+,\times)$ by a $\Delta^0_1$ formula: there is
		a $\Delta^0_1$ formula $\mathrm{Count}(a,b,n,t)$ expressing that $t$ is the
		output of the fixed Turing machine computing $C$ on input $(a,b,n)$ (see
		\cite[Theorem~VI.1.1]{Shoenfield1967}). By construction, this output is
		$|U(a,b)\cap[0,n]|$, which yields the equivalence.
		
		We note that the naive algorithm described above computes $C(a,b,n)$ by
		constructing $U(a,b)$ up to $n$ once and counting membership along the way; in
		particular, $C$ is not only total recursive but primitive recursive, and its
		computation can be bounded by an elementary function of~$n$ for each fixed
		$(a,b)$. We do not attempt to optimise the running time here, as only
		recursiveness is needed for our logical applications.
	\end{proof}
	
	We now analyse the logical form of the natural density inequalities for
	$U(a,b)$. Throughout, we regard the rational parameter $q\in\Q$ as fixed, say
	$q=p/s$ with $p\in\Z$ and $s\in\N$, and we suppress the dependence of
	primitive recursive bounds on this representation.
	
	\begin{theorem}
		\label{thm:density-complexity}
		For each fixed coprime $(a,b)$ and each rational $q\in\Q$, the following
		statements are arithmetical:
		\begin{itemize}[label=\textbullet]
			\item[\emph{(i)}] $\overline{d}(U(a,b))\le q$;
			\item[\emph{(ii)}] $\underline{d}(U(a,b))\ge q$;
			\item[\emph{(iii)}] $d(U(a,b))=q$.
		\end{itemize}
		More precisely, for each fixed $(a,b)$ and $q$, the statements
		\emph{(i)} and \emph{(ii)} can be written as $\Pi^0_3$ sentences of
		arithmetic, and \emph{(iii)} is equivalent to the conjunction of \emph{(i)}
		and \emph{(ii)}, and hence is again $\Pi^0_3$. In particular, all three
		statements are arithmetical and lie (non-optimally) in~$\Delta^0_4$.
	\end{theorem}
	
	\begin{proof}
		Fix coprime $(a,b)$ and rational $q\in\Q$. Let $C(n)=C_{a,b}(n)$ for brevity.
		By Lemma~\ref{lem:counting}, the relation $C(n)=t$ is $\Delta^0_1$, uniformly
		in $n,t$ and the parameters $(a,b)$.
		
		\medskip
		
		\noindent\emph{(i) Upper density $\overline{d}(U(a,b))\le q$.}
		By definition, $\overline{d}(U(a,b))\le q$ means that
		\[
		\limsup_{n\to\infty}\frac{C(n)}{n+1} \le q.
		\]
		Equivalently, for every $k\in\N_{>0}$ there exists $N\in\N$ such that for all
		$n\ge N$ we have
		\[
		\frac{C(n)}{n+1} \le q + \frac{1}{k}.
		\]
		For a fixed rational $q$, the inequality
		\[
		\frac{C(n)}{n+1} \le q + \frac{1}{k}
		\]
		can be rewritten as an inequality between integers using only additions and
		multiplications. Concretely, if $q=p/s$ with $s\ge 1$, then
		\[
		\frac{C(n)}{n+1} \le \frac{p}{s} + \frac{1}{k}
		\quad\Longleftrightarrow\quad
		s k\, C(n) \;\le\; \bigl(p k + s\bigr)\,(n+1),
		\]
		and the right-hand side involves only the parameters $(a,b)$ (through $C(n)$),
		the fixed rational $q$, and the quantified variables $k,n$. Since $C(n)$ is a
		total recursive function, both sides of this inequality are values of
		primitive recursive functions of $(a,b,q,k,n)$, and so the relation
		\[
		\frac{C(n)}{n+1} \le q + \frac{1}{k}
		\]
		is $\Delta^0_1$ as a predicate in $k,n$ (with $(a,b,q)$ as parameters).
		
		Thus the statement ``$\overline{d}(U(a,b))\le q$'' is equivalent to the
		sentence
		\[
		\forall k\,\exists N\,\forall n\ge N\,
		\Phi(k,N,n),
		\]
		where $\Phi(k,N,n)$ is the $\Delta^0_1$ formula expressing the above
		inequality. By the definition of the arithmetical hierarchy, this is a
		$\Pi^0_3$ sentence.
		
		\medskip
		
		\noindent\emph{(ii) Lower density $\underline{d}(U(a,b))\ge q$.}
		By definition, $\underline{d}(U(a,b))\ge q$ means that
		\[
		\liminf_{n\to\infty}\frac{C(n)}{n+1} \ge q.
		\]
		Equivalently, for every $k\in\N_{>0}$ there exists $N\in\N$ such that for all
		$n\ge N$ we have
		\[
		\frac{C(n)}{n+1} \ge q - \frac{1}{k}.
		\]
		As before, for fixed rational $q$ this inequality can be rewritten as an
		inequality between values of primitive recursive functions of $k,n$ (with
		parameters $(a,b,q)$), and hence is expressible by a $\Delta^0_1$ formula
		$\Psi(k,N,n)$ in the language of arithmetic. This gives the sentence
		\[
		\forall k\,\exists N\,\forall n\ge N\,\Psi(k,N,n),
		\]
		which again has quantifier prefix $\forall\exists\forall$ with a $\Delta^0_1$
		matrix. Thus it is a $\Pi^0_3$ sentence.
		
		\medskip
		
		\noindent\emph{(iii) Exact density $d(U(a,b))=q$.}
		By definition, $d(U(a,b))=q$ holds if and only if
		\[
		\overline{d}(U(a,b))\le q
		\quad\text{and}\quad
		\underline{d}(U(a,b))\ge q.
		\]
		By parts (i) and (ii), each of these is equivalent to a $\Pi^0_3$ sentence in
		the language of arithmetic (with parameters $(a,b,q)$). Their conjunction is
		therefore again equivalent to a $\Pi^0_3$ sentence.
		
		In particular, all three statements (i)--(iii) are arithmetical. Since any
		$\Pi^0_3$ sentence can trivially be rewritten as both a $\Sigma^0_4$ and a
		$\Pi^0_4$ sentence by adding dummy quantifiers, each of these statements also
		lies (non-optimally) in $\Delta^0_4$ in the sense that it is equivalent to
		both a $\Sigma^0_4$ and a $\Pi^0_4$ sentence.
	\end{proof}
	
	\begin{figure}
		\centering
		\small
		\setlength{\tabcolsep}{4pt}
		\renewcommand{\arraystretch}{1.1}
		\begin{tabular}{p{0.40\textwidth}p{0.27\textwidth}p{0.27\textwidth}}
			\hline
			Combinatorial statement & For fixed $(a,b)$ & Uniform / family version \\
			\hline
			Finite-window rigidity for $U(1,n)$
			& $\Sigma^0_2$ for each fixed window factor $K$
			& $\Pi^0_3$ for all $K$ (Thm.~\ref{thm:window-rigidity-complexity}) \\[3pt]
			Gap regularity (eventual periodicity)
			& $\Sigma^0_2$ (Thm.~\ref{thm:regularity-sigma02})
			& --- \\[3pt]
			Finitely many even elements
			& $\Sigma^0_2$ (Prop.~\ref{prop:finiteness-evens})
			& --- \\[3pt]
			Upper/lower density \\[1pt]
			$\overline d(U(a,b))\le q$, $\underline d(U(a,b))\ge q$
			& $\Pi^0_3$ (Thm.~\ref{thm:density-complexity})
			& --- \\[3pt]
			Exact density $d(U(a,b))=q$
			& $\Pi^0_3$ (Thm.~\ref{thm:density-complexity})
			& --- \\
			\hline
		\end{tabular}
		\caption{Logical complexity for rigidity, regularity and density statements.}
		\label{fig:rigidity-complexity}
	\end{figure}
	
	% ------------------------------
	% 6. Absoluteness and strength
	% ------------------------------
	
	\section{Forcing absoluteness}
	\label{sec:absoluteness}
	
	\subsection{Absoluteness for arithmetical sentences}
	
	We briefly recall the general principle underlying the absoluteness results we
	obtain for rigidity, regularity and density statements. Let $M$ and $N$ be
	transitive models of $\mathrm{ZFC}$, and write $\N^M$ and $\N^N$ for their
	respective sets of natural numbers. If $\varphi$ is a sentence in the language
	of arithmetic $\mathcal{L}_{\mathrm{PA}}=\{0,1,+,\times,\leq\}$ which is
	\emph{arithmetical} (i.e.\ built using only first-order quantification over
	$\N$ and the arithmetic symbols), then its truth value depends only on the
	structure $(\N,+,\times,\leq)$.
	
	In particular, if $M$ and $N$ have the same natural numbers, i.e.\
	$\N^M=\N^N=\N$, then they agree on the truth of all arithmetical sentences.
	
	This is exactly Lemma~\ref{lem:arithmetical-absoluteness} from
	Section~\ref{sec:preliminaries}, and its proof rests on the fact that the
	satisfaction relation for arithmetical formulas is primitive recursive in the
	underlying structure $(\N,+,\times,\leq)$; see, for instance,
	\cite[Chapter~VI]{Shoenfield1967}; for a related discussion of definability and absoluteness in set theory, 
	see \cite[Chapter~IV and its appendices]{Kunen1980}.
	
	In particular, if $V$ is the ambient universe and $V[G]$ is any forcing
	extension, then $V$ and $V[G]$ have the same natural numbers and the same
	structure $(\N,+,\times,\leq)$, so they agree on all arithmetical sentences.

	\subsection{Absoluteness of rigidity and regularity statements}
	
	We now apply Lemma~\ref{lem:arithmetical-absoluteness} to the specific
	statements analysed in Sections~\ref{sec:rigidity} and
	\ref{sec:regularity-density}. The key input is that all of these statements
	are arithmetical and in fact occupy low levels of the arithmetical hierarchy.
	
	\begin{corollary}
		The following statements are arithmetical and hence absolute between transitive
		models of $\mathrm{ZFC}$ with the same natural numbers:
		\begin{itemize}[label=\textbullet]
			\item for each fixed $K$, the finite-window statement
			$\mathrm{WinRig}(K)$, and the full statement $\mathrm{FWRig}$;
			\item for each fixed coprime $(a,b)$, the regularity statement asserting
			eventual periodicity of the gap sequence of $U(a,b)$;
			\item for each fixed coprime $(a,b)$ and each rational $q\in\Q$, the density
			statements from Theorem~\ref{thm:density-complexity}.
		\end{itemize}
	\end{corollary}
	
	\begin{proof}
		Theorem~\ref{thm:window-rigidity-complexity} shows that
		$\mathrm{WinRig}(K)$ is $\Sigma^0_2$ for fixed $K$ and that
		$\mathrm{FWRig}$ is $\Pi^0_3$. Theorem~\ref{thm:regularity-sigma02}
		shows that regularity of $U(a,b)$ is equivalent to a $\Sigma^0_2$ sentence.
		Finally, Theorem~\ref{thm:density-complexity} shows that the upper and lower
		density inequalities, and hence the exact density statements, are equivalent
		to $\Pi^0_3$ sentences for fixed $(a,b)$ and rational $q$.
		
		In particular, in each case the statement in question is arithmetical (i.e.\
		$\Sigma^0_n$ or $\Pi^0_n$ for some $n$). Lemma
		\ref{lem:arithmetical-absoluteness} therefore implies that any two
		transitive models of $\mathrm{ZFC}$ with the same natural numbers agree on the
		truth value of each of these statements.
	\end{proof}
	
	\begin{remark}
		The transitivity assumption in the preceding corollary can be relaxed. If $M$
		and $N$ are $\omega$-models of $\mathrm{ZFC}$ (i.e.\ their natural numbers are
		isomorphic to the standard $\N$), then they agree on all arithmetical
		sentences, regardless of transitivity. In particular, the rigidity, regularity
		and density statements considered in this paper have the same truth value in
		any two $\omega$-models of $\mathrm{ZFC}$, and hence are insensitive to
		variations in the higher set-theoretic universe even in that more general
		sense.
	\end{remark}
	
	\begin{corollary}
		Set forcing cannot change the truth value of any of the preceding
		arithmetical statements.
	\end{corollary}
	
	\begin{proof}
		Let $\varphi$ be any of the statements listed in the preceding corollary.
		Suppose $V$ is a model of $\mathrm{ZFC}$ and $V[G]$ is a forcing extension
		obtained by set forcing. Then $V$ and $V[G]$ are both transitive models of
		$\mathrm{ZFC}$ with the same natural numbers: forcing does not create new
		finite ordinals, so $\omega^{V} = \omega^{V[G]}$ and $(\N,+,\times,\leq)$ is
		the same structure in both models. By the previous corollary, $V$ and $V[G]$
		must agree on the truth value of~$\varphi$.  Equivalently, if $V\models
		\varphi$ then the top condition forces $\varphi$, and if $V\models\neg
		\varphi$ then the top condition forces $\neg\varphi$.
	\end{proof}

	\begin{remark}[Truth versus provability]
		Forcing invariance is a semantic statement about models with the same
		standard natural numbers.  It does not imply that extensions of set theory
		by large-cardinal axioms are arithmetically conservative, and it does not
		preclude set-theoretic arguments from proving arithmetic theorems.  No such
		proof-theoretic claim is made here.
	\end{remark}

	\subsection{Formalisation over weak arithmetic}

	The displayed complexity bounds are syntactic upper bounds in the standard
	model.  A reverse-mathematical calibration would additionally have to verify,
	inside the chosen base theory, the totality of the particular Ulam recursion
	and the equivalence between the coded predicates and the intended sequences.
	We do not claim here that every preceding proof formalises ``without change''
	in $\mathsf{RCA}_0$; that is a separate question.

	% ------------------------------
	% 7. Model-theoretic consequences
	% ------------------------------
	
	\section{Model-theoretic aspects of \texorpdfstring{$(\N,+,\Ulam_{a,b})$}{(N,+,Ulam)}}
	\label{sec:model-theory}
	
	\subsection{Model-theoretic background}
	\label{sec:model-theory-background}
	
	We briefly recall the model-theoretic notions that appear in
	Section~\ref{sec:model-theory}. We work in the usual framework of first-order
	model theory; see, for example, standard texts on stability and NIP theories
	for further background.
	
	\begin{definition}[Independence property (IP) and NIP]
		Let $T$ be a complete first-order theory and $\varphi(x;y)$ a formula in the
		language of $T$, where $x$ and $y$ are tuples of variables. We say that
		$\varphi(x;y)$ has the \emph{independence property (IP)} if there exists a
		model $M\models T$ and a sequence of parameters
		$(b_i)_{i\in\N}$ from $M^{|y|}$ such that for every subset
		$S\subseteq\N$ there is $a_S\in M^{|x|}$ with
		\[
		M \models \varphi(a_S;b_i) \quad\Longleftrightarrow\quad i\in S
		\qquad\text{for all } i\in\N.
		\]
		In words, the family of definable sets
		$\{\varphi(M;b_i) : i\in\N\}$ realises all patterns on $\N$.
		
		A theory $T$ is said to be \emph{NIP} (``not the independence property'') if
		no formula $\varphi(x;y)$ has IP in any model of $T$. Equivalently, every
		$\varphi(x;y)$ has \emph{finite} Vapnik--Chervonenkis dimension in the sense
		recalled below.
	\end{definition}
	
	\begin{definition}[VC-dimension]
		Let $X$ be a set and $\mathcal{F}\subseteq\mathcal{P}(X)$ a family of subsets
		of $X$. A finite subset $A\subseteq X$ is said to be \emph{shattered} by
		$\mathcal{F}$ if
		\[
		\{ F\cap A : F\in\mathcal{F} \} = \mathcal{P}(A),
		\]
		i.e.\ every subset of $A$ arises as the intersection of $A$ with some
		$F\in\mathcal{F}$. The \emph{VC-dimension} of $\mathcal{F}$, written
		$\mathrm{VCdim}(\mathcal{F})$, is the supremum of the cardinalities $|A|$ of
		finite sets $A\subseteq X$ that are shattered by~$\mathcal{F}$ (or $\infty$ if
		no such finite bound exists).
		
		Given a formula $\varphi(x;y)$ and a structure $M$, we consider the family
		\[
		\mathcal{F}_\varphi = \{ \varphi(M;b) : b\in M^{|y|} \}
		\]
		of subsets of $M^{|x|}$. We say that $\varphi$ has finite VC-dimension in $M$
		if $\mathrm{VCdim}(\mathcal{F}_\varphi)<\infty$, and we say that $T$ is NIP
		iff every formula has finite VC-dimension in every model of $T$.
	\end{definition}
	
	\begin{definition}[dp-minimality]
		A complete NIP theory $T$ is called \emph{dp-minimal} if, informally, every
		definable one-parameter family behaves like a ``linear'' family from the point
		of view of dependence. There are several equivalent formal definitions; we
		sketch one that suffices for our purposes.
		
		We say that $T$ is \emph{dp-minimal} if for every model $M\models T$ and every
		formula $\varphi(x;y)$ with $|x|=1$, one cannot find two sequences of
		parameters $(b_i^0)_{i\in\N}$ and $(b_j^1)_{j\in\N}$ from $M^{|y|}$ and
		elements $(a_{i,j})_{i,j\in\N}$ of $M$ such that
		\[
		M \models \varphi(a_{i,j};b_k^\ell)
		\quad\Longleftrightarrow\quad
		\bigl( \ell=0\ \wedge\ k=i \bigr) \text{ or } \bigl( \ell=1\ \wedge\ k=j \bigr)
		\]
		for all $i,j,k\in\N$ and $\ell\in\{0,1\}$. In other words, there is no
		``independent pattern'' of depth~$2$ in a single variable. For more systematic
		treatments and equivalent characterisations (via dp-rank, VC-density, etc.),
		one may consult the standard literature on NIP and dp-minimal theories.
		
		A basic example is Presburger arithmetic: the complete theory
		$\Th(\N,+,0,1)$ of addition on the natural numbers is NIP and dp-minimal,
		and its definable sets are semilinear.
	\end{definition}
	
	\begin{definition}[Interpretation and definitional expansion]
		Let $\mathcal{M}$ and $\mathcal{N}$ be first-order structures. We say that
		$\mathcal{N}$ is \emph{interpretable} in $\mathcal{M}$ if there is a
		definable set $D\subseteq M^k$ (for some $k$) together with definable
		relations and functions on $D$ making $(D,\dots)$ isomorphic to $\mathcal{N}$.
		We say that a theory $T$ \emph{interprets} a structure $\mathcal{N}$ if some
		(or equivalently every) model of $T$ does.
		
		An expansion $\mathcal{M}'$ of a structure $\mathcal{M}$ is called a
		\emph{definitional expansion} if every new symbol in the language of
		$\mathcal{M}'$ is definable in $\mathcal{M}$, and conversely every
		$\mathcal{M}'$-formula is equivalent (in $\Th(\mathcal{M})$) to a formula in
		the original language. Definitional expansions preserve model-theoretic
		properties such as NIP, dp-minimality, and interpretability power.
	\end{definition}
	
	For general background on NIP theories, VC-dimension, and dp-minimality we refer to Simon’s monograph \cite{SimonNIP}.
	
	\subsection{The structure \texorpdfstring{$(\N,+,\Ulam_{a,b})$}{(N,+,Ulam)}}
	
	For each fixed coprime pair $(a,b)$ we consider the expansion of Presburger
	arithmetic by a unary predicate naming the Ulam sequence $U(a,b)$.
	
	\begin{definition}
		For fixed coprime $(a,b)$ we write
		\[
		\Str{U}_{a,b} = (\N,+,0,1,\Ulam_{a,b}),
		\]
		where $\Ulam_{a,b}(m)$ abbreviates $\Ulam(a,b,m)$, the $\Delta^0_1$ predicate
		from Theorem~\ref{thm:uniform-ulam} defining membership in $U(a,b)$.
	\end{definition}
	
	The reduct of $\Str{U}_{a,b}$ to the language $\{0,1,+\}$ is the standard
	model of Presburger arithmetic, which admits quantifier elimination and whose
	definable subsets of $\N^k$ are exactly the \emph{semilinear} sets, i.e.\
	finite unions of linear sets and singletons. In particular, every definable
	subset of $\N$ is a finite union of arithmetic progressions and isolated
	points. The expansion by $\Ulam_{a,b}$ need not preserve this semilinear
	behaviour in general: for an arbitrary recursive subset $P\subseteq\N$, the
	structure $(\N,+,0,1,P)$ may define sets far more complicated than semilinear
	ones, and may even interpret multiplication.

	For the model-theoretic conclusions we use only a condition on one fixed
	sequence.

	\begin{definition}[Structural conditions]
		\label{def:rigidity-hierarchy}
		Let $(a,b)$ be coprime with $a<b$.  We consider:
		\begin{enumerate}[label=\emph{(R\arabic*)}]
			\item \emph{Gap regularity}: there are $N$ and $p\geq1$ such that
			$g_{k+p}=g_k$ for every $k\geq N$.
			\item \emph{Bounded gaps}: there is $B$ such that $g_k\leq B$ for all $k$.
			\item \emph{Density existence}: $d(U(a,b))$ exists.
			\item \emph{Positive lower density}: $\underline d(U(a,b))>0$.
		\end{enumerate}
	\end{definition}

	Gap regularity implies bounded gaps, existence of density, and positive lower
	density.  Bounded gaps imply positive lower density.  No implication from
	finite-window rigidity of the varying family $U(1,n)$ to gap regularity of a
	fixed $U(a,b)$ is asserted.

	We now establish the key consequences of gap regularity.
	
	\begin{proposition}[Consequences of gap regularity]
		\label{prop:gap-regularity-consequences}
		Suppose $U(a,b)$ satisfies condition \emph{(R1)}: the gap sequence $(g_k)$ is
		eventually periodic with period $p$ and threshold $N$. Then:
		\begin{enumerate}[label=\emph{(\roman*)}]
			\item The density $d(U(a,b))$ exists and equals
			\[
			d(U(a,b)) = \frac{p}{\sum_{j=0}^{p-1}g_{N+j}}.
			\]
			\item The set $U(a,b)$ is the union of a finite set and a set that is
			Presburger-definable with parameters.
			\item The expansion $\Str{U}_{a,b}$ is NIP.
		\end{enumerate}
	\end{proposition}
	
	\begin{proof}
		\emph{(i)} For $k\ge N$, write $k = N + mp + r$ with $0 \le r < p$ and
		$m \ge 0$. The partial sums of gaps satisfy
		\[
		u_{N+mp+r} = u_N + m \cdot G + \sum_{j=0}^{r-1} g_{N+j},
		\]
		where $G = \sum_{j=0}^{p-1} g_{N+j}$ is the sum of one period of gaps. For
		large $n$, the number of Ulam elements in $[u_N, n]$ is asymptotically
		$p \cdot (n - u_N)/G$. It follows that
		\[
		\lim_{n\to\infty} \frac{|U(a,b)\cap[0,n]|}{n} = \frac{p}{G}.
		\]
		
		\emph{(ii)} The tail $\{u_k : k \ge N\}$ consists of elements of the form
		\[
		u_N + mG + \sum_{j=0}^{r-1} g_{N+j}
		\]
		for $m \ge 0$ and $0 \le r < p$. For each fixed $r$, this is an arithmetic
		progression with first term $u_N + \sum_{j=0}^{r-1} g_{N+j}$ and common
		difference $G$. Thus the tail is a union of $p$ arithmetic progressions.
		Arithmetic progressions are Presburger-definable, so the tail is
		Presburger-definable with parameters $(u_N, g_N, \dots, g_{N+p-1})$. The
		finite initial segment $\{u_k : k < N\}$ is also Presburger-definable (as
		a finite set of constants).
		
		(iii) By (ii), there is a formula $\psi(x,\bar c)$ in the language $\{0,1,+\}$
		with a finite tuple of parameters $\bar c$ such that $m\in U(a,b)$ iff
		$(\N,+,0,1)\models\psi(m,\bar c)$. Naming the parameters $\bar c$ gives a
		(definitional) expansion of Presburger arithmetic that remains NIP (since NIP is
		preserved under naming constants). Now adding a new unary predicate symbol
		$P(x)$ together with the defining axiom $\forall x\,(P(x)\leftrightarrow
		\psi(x,\bar c))$ is again a definitional expansion, so the resulting theory is
		NIP. Interpreting $P$ as $\Ulam_{a,b}$ yields $\Str U_{a,b}$, which is therefore
		NIP.
	\end{proof}
	
	\begin{remark}[Gap regularity and finite-window rigidity]
		\label{rem:gap-vs-strong}
		Finite-window rigidity concerns uniform descriptions of the varying family
		$U(1,n)$ on intervals $[1,Kn]$; gap regularity concerns the tail of one fixed
		sequence.  Neither statement should be substituted for the other.  Under gap
		regularity, the displayed formula in
		Proposition~\ref{prop:gap-regularity-consequences} shows directly that the
		tail is a finite union of arithmetic progressions.
	\end{remark}
	
	We now record the main model-theoretic consequence. Recall that definitional 
	expansions do not change the class of definable sets (up to naming additional 
	definable predicates) or the interpretability power of a theory: any structure
	interpretable in a definitional expansion is already interpretable in the
	original structure.
	
	\begin{proposition}[Model-theoretic tameness under gap regularity]
		\label{prop:no-multiplication}
		Assume that $U(a,b)$ satisfies gap regularity~\textup{(R1)} of
		Definition~\ref{def:rigidity-hierarchy}, i.e., its gap sequence $(g_k)$
		is eventually periodic. Then the theory $\Th(\Str{U}_{a,b})$ is a
		definitional expansion of Presburger arithmetic. In particular:
		\begin{itemize}[label=\textbullet]
			\item $\Th(\Str{U}_{a,b})$ does not interpret the full arithmetic structure
			$(\N,+,\times)$;
			\item the expansion $\Str{U}_{a,b}$ is NIP and dp-minimal; in particular,
			each fixed definable family of subsets of $\N$ has finite VC-dimension.
		\end{itemize}
	\end{proposition}
	
	\begin{proof}
		By Proposition~\ref{prop:gap-regularity-consequences}\emph{(ii)}, under
		the assumption of gap regularity~\textup{(R1)} there is a formula
		$\varphi(x)$ in the language $\{0,1,+\}$ with parameters such that for
		all $m\in\N$,
		\[
		\Str{U}_{a,b}\models \Ulam_{a,b}(m) \quad\Longleftrightarrow\quad
		(\N,+,0,1)\models \varphi(m).
		\]
		In other words, $U(a,b)$ is definable (up to a finite initial segment,
		which is also definable) in the underlying Presburger structure. Thus
		the predicate symbol $\Ulam_{a,b}$ is definable in the reduct
		$(\N,+,0,1)$, and $\Str{U}_{a,b}$ is a definitional expansion of
		$(\N,+,0,1)$: every formula in the expanded language can be translated
		into an equivalent formula in the original Presburger language by
		replacing atomic instances of $\Ulam_{a,b}(x)$ with~$\varphi(x)$.
		
		Let $T_{\mathrm{Pres}}$ denote the complete theory of $(\N,+,0,1)$ and
		let $T_{U}$ denote $\Th(\Str{U}_{a,b})$. Since $T_{U}$ is a
		conservative definitional extension of $T_{\mathrm{Pres}}$, any
		structure interpretable in $T_{U}$ is already interpretable in
		$T_{\mathrm{Pres}}$.
		
		It is a classical fact that Presburger arithmetic does not interpret the
		full ring $(\N,+,\times)$; in particular, multiplication is not
		interpretable in $(\N,+,0,1)$. If $(\N,+,\times)$ were interpretable in
		$T_{U}$, then by definitionality it would be interpretable in
		$T_{\mathrm{Pres}}$, a contradiction. This proves that $T_{U}$ does not
		interpret $(\N,+,\times)$.
		
		Finally, Presburger arithmetic is NIP and dp-minimal, and NIP and
		dp-minimality are preserved under definitional expansions. Hence
		$\Str{U}_{a,b}$ is again NIP and dp-minimal, and, in particular, every
		family of definable subsets of $\N$ in $\Str{U}_{a,b}$ has uniformly
		bounded VC-dimension. Since Presburger arithmetic is decidable, the
		same holds for $\Th(\Str{U}_{a,b})$.
	\end{proof}
	
	\medskip
	
	In Section~\ref{subsec:unconditional} we ask how far one can push such
	tameness without assuming gap regularity.
	
	\begin{remark}[Computational consequences under rigidity]
		\label{rem:computational}
		Under the hypotheses of Proposition~\ref{prop:no-multiplication}, every
		first-order sentence in the language of $\Str{U}_{a,b}$ can be effectively
		translated into an equivalent sentence of Presburger arithmetic, uniformly in
		the parameters $(a,b)$ and in the pattern code witnessing rigidity. In
		particular, the theory $\Th(\Str{U}_{a,b})$ is decidable in this regime, and
		one may in principle combine standard decision procedures for Presburger
		arithmetic (such as Cooper's algorithm or automata-theoretic methods) with the
		arithmetisation of $U(a,b)$ to obtain explicit algorithms for answering
		first-order questions about $U(a,b)$. We do not attempt to optimise such
		algorithms here, as our primary focus is on logical complexity rather than
		computational complexity.
	\end{remark}
	
	There is a substantial body of work on decidable and undecidable expansions of Presburger arithmetic; see, for example, Point's study of decidable extensions \cite{Point2000DecidablePresburger} and B\`es' survey on arithmetical definability \cite{Bes2001SurveyArithDef}.
	
	\subsection{Unconditional and partial results}
	\label{subsec:unconditional}
	
	The preceding results establish model-theoretic tameness of $\Str{U}_{a,b}$
	under the hypothesis of gap regularity~\textup{(R1)}. In this subsection we investigate what can be
	said without these assumptions, either unconditionally or under weaker
	structural hypotheses.
	While full tameness results remain out of reach
	in the general case, we identify several constraints that the Ulam construction
	imposes on the expansion $\Str{U}_{a,b}$ and formulate precise conjectures
	about intermediate cases.
	
	\subsubsection{The unique-sum property and its consequences}
	
	The defining feature of Ulam sequences---that each term beyond the initial
	pair has a unique representation as a sum of two distinct earlier terms---is
	a strong structural constraint. We begin by formalising this and extracting
	some unconditional consequences.
	
	\begin{definition}[Representation function]
		\label{def:representation-function}
		For a set $P\subseteq\N$ and $n\in\N$, the \emph{representation function}
		$r_P(n)$ counts the number of ways to write $n$ as a sum of two distinct
		elements of $P$:
		\[
		r_P(n) = \bigl|\{(x,y)\in P^2 : x<y \text{ and } x+y=n\}\bigr|.
		\]
	\end{definition}
	
	\begin{lemma}[Representation bound for Ulam sequences]
		\label{lem:representation-bound}
		Let $(a,b)$ be coprime with $a<b$. For every $m\in U(a,b)\setminus\{a,b\}$
		we have $r_{U(a,b)}(m)=1$. Moreover, for every $n\notin U(a,b)$ with $n>b$,
		either $r_{U(a,b)\cap[1,n)}(n)=0$ or $r_{U(a,b)\cap[1,n)}(n)\ge 2$.
	\end{lemma}
	
	\begin{proof}
		The first claim is immediate from the definition of the Ulam sequence: an
		integer $m>b$ is adjoined to $U(a,b)$ at stage $k$ if and only if $m$ is
		the smallest integer exceeding all previously included terms and having
		exactly one representation as a sum of two distinct earlier terms. Thus
		$r_{U(a,b)\cap[1,m)}(m)=1$ at the moment of inclusion, and since no further
		terms less than $m$ are added subsequently, this remains $r_{U(a,b)}(m)=1$
		in the final sequence.
		
		For the second claim, suppose $n>b$ and $n\notin U(a,b)$. At the stage when
		the Ulam construction considers $n$, either $n$ has no representations as a
		sum of two distinct earlier terms (so $r_{U(a,b)\cap[1,n)}(n)=0$), or $n$
		has at least two such representations (so $r_{U(a,b)\cap[1,n)}(n)\ge 2$). In
		the latter case, since the Ulam sequence $U(a,b)\cap[1,n)$ is already
		determined by stage $n$, this representation count persists.
	\end{proof}
	
	The representation bound suggests heuristic constraints on the density of
	Ulam sequences, though rigorous upper bounds appear difficult to establish.
	
	\begin{remark}[Density heuristics]
		\label{rem:density-heuristics}
		Intuitively, one expects the density of any Ulam sequence to be bounded
		well away from $1$. The reason is as follows: if $U = U(a,b)$ had density
		$d$ close to $1$, then the number of pairs from $U \cap [1,n]$ would be
		$\binom{dn}{2} \approx d^2 n^2 / 2$, while the range of possible sums is
		only $[a+b, 2n]$, an interval of length $< 2n$. By the pigeonhole principle,
		many integers would have multiple representations, and by
		Lemma~\ref{lem:representation-bound}, such integers cannot belong to $U$
		(beyond the initial terms). This suggests a tension between high density
		and the unique-representation requirement.
		
		Making this precise is nontrivial. A naive application of the Cauchy--Schwarz
		inequality to the second moment $\sum_m r_{U_n}(m)^2$ yields only weak bounds
		on the number of multiply-represented integers, insufficient to establish
		a sharp density constraint such as $\overline{d}(U) \le 1/2$. The difficulty
		lies in controlling the distribution of representation counts among
		non-members: while each must have $0$ or $\ge 2$ representations, the
		total ``energy'' absorbed by non-members depends on how representation
		counts are distributed, and crude bounds lose too much information.
		
		Numerical evidence suggests that densities of Ulam sequences, when they
		exist, are quite small. For the classical sequence $U(1,2)$, computations
		indicate $d(U(1,2)) \approx 0.07$ \cite{Steinerberger2017HiddenSignal},
		and other well-studied pairs yield similarly small values. The gap between
		these empirical densities and even a conjectural bound of $1/2$ suggests
		that the rigidity phenomena impose much stronger constraints than naive
		counting arguments capture.
	\end{remark}
	
	\subsubsection{The additive character of the Ulam construction}
	
	We now turn to structural features of the expansion $\Str{U}_{a,b}$ that
	follow from the purely additive nature of the Ulam construction.
	
	\begin{definition}[Additively closed construction]
		\label{def:additively-closed}
		A subset $P\subseteq\N$ is said to arise from an \emph{additively closed
			construction} if there is a recursive procedure defining $P$ in which
		membership of $n$ in $P$ depends only on:
		\begin{itemize}
			\item the values of $P\cap[1,n)$ (the initial segment);
			\item linear combinations and comparisons of elements of $P\cap[1,n)$;
			\item bounded quantification over $P\cap[1,n)$.
		\end{itemize}
		In particular, the procedure never references multiplication, division,
		exponentiation, or other non-additive operations on the elements themselves.
	\end{definition}
	
	\begin{lemma}
		\label{lem:ulam-additively-closed}
		For every coprime $(a,b)$, the Ulam sequence $U(a,b)$ arises from an
		additively closed construction in the sense of
		Definition~\ref{def:additively-closed}.
	\end{lemma}
	
	\begin{proof}
		The Ulam construction proceeds as follows: given $U(a,b)\cap[1,n)$, an
		integer $n$ is adjoined if and only if (i) $n>\max(U(a,b)\cap[1,n))$, and
		(ii) there is exactly one pair $(x,y)$ with $x,y\in U(a,b)\cap[1,n)$,
		$x<y$, and $x+y=n$. Condition (i) involves only comparison. Condition (ii)
		involves only addition of elements and counting (bounded quantification).
		No multiplicative operations appear. Thus the construction is additively
		closed.
	\end{proof}
	
	The significance of this observation is that additively closed constructions
	produce predicates whose interaction with the additive structure of $\N$ is
	constrained. While we cannot prove unconditionally that $\Str{U}_{a,b}$ fails
	to interpret multiplication, we can articulate why such an interpretation
	would be surprising.
	
	\begin{proposition}[Decidability of the atomic diagram]
		\label{prop:atomic-decidable}
		Let $(a,b)$ be coprime with $a<b$. The atomic diagram of
		$\Str{U}_{a,b}=(\N,+,0,1,\Ulam_{a,b})$ is decidable: given any
		quantifier-free sentence $\varphi$ with parameters from $\N$, one can
		effectively determine whether $\Str{U}_{a,b}\models\varphi$.
	\end{proposition}
	
	\begin{proof}
		The predicate $\Ulam(a,b,m)$ is $\Delta^0_1$ by Theorem~\ref{thm:uniform-ulam},
		hence decidable. Any quantifier-free formula in the language
		$\{0,1,+,\Ulam_{a,b}\}$ with parameters involves only finitely many membership
		queries to $U(a,b)$ and finitely many linear comparisons, all of which are
		decidable.
	\end{proof}
	
	\begin{remark}[On Presburger definability]
		\label{rem:presburger-definability}
		Presburger-definable subsets of $\N$ are exactly the eventually periodic sets:
		finite unions of arithmetic progressions and singletons. It is expected, in
		line with computational evidence and the conjectured lack of eventual
		periodicity of gaps for most pairs $(a,b)$, that $\Ulam_{a,b}$ is not
		definable in $(\N,+,0,1)$ unless its characteristic set is eventually periodic.
		Indeed, definability of $U(a,b)$ in Presburger arithmetic would force eventual
		periodicity of its characteristic function, which appears to fail for generic
		pairs based on extensive numerical data.
	\end{remark}
	
	\begin{remark}[On quantifier elimination]
		\label{rem:quantifier-elimination}
		It seems unlikely that $\Str{U}_{a,b}$ admits quantifier elimination in the
		language $\{0,1,+,\Ulam_{a,b}\}$ when $U(a,b)$ is not eventually periodic.
		Already the set $\{(x,y)\in\N^2 : x+y\in U(a,b)\}$ is definable by a simple
		existential formula, but describing it quantifier-free would require the
		predicate $\Ulam_{a,b}$ to interact with addition in a highly structured way.
		We do not pursue a definitive result on quantifier elimination here.
	\end{remark}
	
	\subsubsection{Comparison with known interpretations of arithmetic}
	
	To place the model-theoretic status of $\Str{U}_{a,b}$ in context, we briefly
	recall when an expansion $(\N,+,P)$ interprets full arithmetic.
	
	\begin{definition}[Multiplicative encoding]
		\label{def:multiplicative-encoding}
		We say that a predicate $P\subseteq\N$ \emph{encodes multiplication} if
		there exist first-order formulas $\varphi_\times(x,y,z)$ and $\varphi_D(x)$
		in the language $\{0,1,+,P\}$ such that $\varphi_D$ defines an infinite
		subset $D\subseteq\N$ and
		\[
		(\N,+,P)\models \varphi_\times(a,b,c) \quad\Longleftrightarrow\quad
		a\cdot b = c
		\]
		for all $a,b,c\in D$. Equivalently, $(\N,+,P)$ interprets a nontrivial
		fragment of $(\N,+,\times)$.
	\end{definition}
	
	Predicates defined in terms of strong multiplicative or exponential growth
	are known to produce undecidable expansions where arithmetic is interpretable.
	Classical examples include:
	\begin{itemize}
		\item $P = \{k^n : n\in\N\}$ for a fixed base $k\ge 2$ (powers);
		\item $P = \{n! : n\in\N\}$ (factorials);
		\item $P = \{p_n : n\in\N\}$ (the sequence of primes).
	\end{itemize}
	For these and related results, see Villemaire~\cite{Villemaire1992},
	Muchnik--Semenov~\cite{MuchnikSemenov2003}, Point~\cite{Point2000DecidablePresburger},
	and the survey by B\`es~\cite{Bes2001SurveyArithDef}. The common feature of
	these examples is that the predicate $P$ encodes multiplicative information
	through its growth pattern or divisibility structure.

	Beyond the positive-period gap-regular case, no NIP or non-interpretability
	claim is proved here.  Determining whether arbitrary Ulam expansions interpret
	full arithmetic remains an open model-theoretic problem.

	\subsubsection{Effectivity and computability considerations}
	
	We conclude this subsection with remarks on the effective content of the
	results under partial rigidity. The following is essentially a repackaging
	of the gap-regularity analysis in a more computational form.
	
	\begin{proposition}[Effective consequences of given gap witnesses]
		\label{prop:effective-bounds}
		Suppose a threshold $N$ and a positive period $p\geq1$ are given such that
		$g_{k+p}=g_k$ for all $k\ge N$. Then:
		\begin{enumerate}[label=\emph{(\roman*)}]
			\item The density $d(U(a,b))$ is computable (as a rational number) from the
			witnesses.
			\item The theory $\Th(\Str{U}_{a,b})$ is decidable.
			\item Membership in any definable subset of $\N^k$ is decidable, uniformly
			in the defining formula.
		\end{enumerate}
	\end{proposition}
	
	\begin{proof}
		\emph{(i)} From $N$ and $p$, compute $G = \sum_{j=0}^{p-1} g_{N+j}$ by
		simulating the Ulam construction up to stage $N+p$. Then
		$d(U(a,b)) = p/G$ by Proposition~\ref{prop:gap-regularity-consequences}(i).
		
		\emph{(ii)} By Proposition~\ref{prop:gap-regularity-consequences}(ii), the
		predicate $\Ulam_{a,b}$ is Presburger-definable with parameters determined
		by the witnesses. A formula in the language of $\Str{U}_{a,b}$ can therefore
		be translated into an equivalent Presburger formula (with the parameter
		values substituted). Presburger arithmetic is decidable (e.g.\ via Cooper's
		algorithm or automata-theoretic methods), so the translated formula can be
		decided.
		
		\emph{(iii)} Given a formula $\varphi(\bar x)$ in the language of
		$\Str{U}_{a,b}$ and a tuple $\bar n\in\N^k$, we can decide
		$\Str{U}_{a,b}\models\varphi(\bar n)$ by the procedure in~(ii).
	\end{proof}
	
	\begin{remark}
		For any one fixed gap-regular pair, particular witnesses may be hard-coded,
		so the preceding consequences are effective non-uniformly.  A distinct
		question is whether one algorithm can recover witnesses from $(a,b)$ uniformly
		over a specified class of pairs.  Gap regularity by itself does not supply
		such a uniform witness-extraction procedure.
	\end{remark}
	
	% ------------------------------
	% 8. Further directions
	% ------------------------------
	
	\section{Further directions and open problems}
	\label{sec:further-directions}
	
	The results of this paper place a range of natural statements about Ulam
	sequences $U(a,b)$ at low levels of the arithmetical hierarchy, and use this
	to derive absoluteness and basic model-theoretic consequences for the
	expansions $(\N,+,\Ulam_{a,b})$. We end by listing some directions in which
	this picture could be refined or extended.
	
	\subsection*{Refining complexity bounds}
	
	For rigidity, regularity and density, we have shown low arithmetical upper
	bounds: finite-window rigidity is $\Sigma^0_2$ for a fixed window factor and
	$\Pi^0_3$ uniformly over all factors; positive-period gap regularity is
	$\Sigma^0_2$ for fixed $(a,b)$; and the density assertions considered here
	are $\Pi^0_3$.
	Several natural questions remain.
	
	\begin{itemize}[label=\textbullet]
		\item Are the corresponding \emph{index sets} of parameter pairs or pattern
		codes complete at these levels under many-one reductions?  Completeness is a
		property of such sets of instances, not of one fixed truth-valued sentence.
		\item Can one show that certain weaker forms of rigidity, or certain
		regularity properties, are in fact equivalent to simpler
		$\Sigma^0_1$ or $\Pi^0_1$ sentences, perhaps after appropriate
		reformulation?
		\item Conversely, are there natural Ulam-like constructions for which
		combinatorial conjectures give rise to statements of genuinely higher
		arithmetical complexity?
	\end{itemize}
	
	A systematic proof-theoretic analysis by reverse-mathematical methods would also be of interest,
	for instance by calibrating $\mathrm{FWRig}$ or
	uniform regularity principles between $\mathsf{RCA}_0$ and stronger systems.
	
	\subsection*{Other Ulam-like constructions}
	
	The framework developed here applies to any additive construction whose
	membership relation is recursive and whose conjectured large-scale structure
	can be encoded by finite patterns with linear dependence on parameters.
	Natural candidates for similar analysis include:
	
	\begin{itemize}[label=\textbullet]
		\item \emph{Higher-dimensional Ulam sets} in $\Z^n$, as studied by 
		Kravitz--Steinerberger~\cite{KravitzSteinerberger2018} and 
		Hinman--Kuca--Schlesinger--Sheydvasser~\cite{HinmanKucaSchlesingerSheydvasser2019Involve}. 
		The arithmetisation methods extend naturally to recursively presented 
		ordered abelian groups.
		\item Variants in which the ``unique representation'' condition is replaced
		by bounded multiplicity (e.g.\ at most $k$ representations) or by
		weighted counting of representations of different types.
		\item Other additive sequences such as sum-free sets, Sidon sets, and more
		general $B_h$-sequences, where membership is again defined via uniqueness
		or bounded multiplicity of representations as sums of previous terms.
		Many of these admit conjectural or proven interval-with-mask descriptions,
		and it would be natural to ask whether their rigidity and density
		statements enjoy analogous arithmetical and absoluteness properties.
		\item Polynomial Ulam sequences in the sense of Sheydvasser
		\cite{Sheydvasser2021LinearPolynomials}, and more general algebraic
		deformations where the parameters $(a,b)$ are replaced by integer
		polynomials or higher-dimensional algebraic data.
	\end{itemize}
	
	For such systems, one can ask to what extent the sort of arithmetisation and
	complexity analysis developed here still applies, and whether similar
	absoluteness phenomena persist.
	
	\subsection*{Model-theoretic classification}
	
	The expansion $(\N,+,\Ulam_{a,b})$ provides a concrete test case for the
	model theory of expansions of Presburger arithmetic by unary predicates. Under
	positive-period gap regularity, we showed that $\Str{U}_{a,b}$ is a definitional expansion of
	Presburger arithmetic and hence is NIP and dp-minimal
	(Proposition~\ref{prop:no-multiplication}). Several directions suggest
	themselves.
	
	\begin{itemize}[label=\textbullet]
		\item Can one prove NIP or other tameness properties of $\Str{U}_{a,b}$ under
		substantially weaker combinatorial hypotheses than eventual periodicity of
		the gap sequence?
		\begin{conjecture}[Tameness of all Ulam expansions]
			\label{conj:tame-ulam}
			For every coprime $(a,b)$, the theory $\Th(\Str{U}_{a,b})$ is NIP and does not
			interpret $(\N,+,\times)$. In particular, all expansions $(\N,+,\Ulam_{a,b})$
			are model-theoretically tame in the sense of classification theory.
		\end{conjecture}
		\item At the family level, can one obtain uniform bounds on VC-dimension or
		dp-rank for classes of sets definable in the multi-sorted structure
		obtained by varying $(a,b)$ and retaining a uniform predicate
		$\Ulam(a,b,m)$?
	\end{itemize}
	
	More broadly, it would be interesting to relate the behaviour of Ulam
	sequences to existing general results on expansions of Presburger arithmetic
	by sparse sets or sets of prescribed asymptotic density.
	
	\begin{remark}[Automatic structures]
		\label{rem:automatic}
		Expansions of $(\N,+)$ whose definable sets are semilinear are closely related
		to \emph{automatic structures}, where elements and relations are recognised by
		finite automata. Under positive-period gap regularity, the predicate $\Ulam_{a,b}$ is
		Presburger-definable, so $(\N,+,\Ulam_{a,b})$ sits very near this automatic
		realm. It would be interesting to make this connection precise, for example by
		asking whether suitable presentations of $\Str{U}_{a,b}$ are automatic in the
		sense of Khoussainov--Nerode, and whether the pattern decompositions arising
		from Ulam rigidity can be recovered by purely automata-theoretic means.
	\end{remark}
	
	\subsection*{Broader logical context}
	
	Finally, the arithmetical nature of the Ulam problems suggests possible
	connections with other areas of logic. For example:
	
	\begin{itemize}[label=\textbullet]
		\item Since all the statements considered here are arithmetical, their truth values
		are invariant under set forcing and do not depend on the ambient set-theoretic
		universe. Thus any genuine independence phenomena would have to be of a
		proof-theoretic nature (for example, independence from a weak theory of
		arithmetic such as $\mathsf{RCA}_0$ or $\mathsf{PA}$), which are reflected in
		the behaviour of these statements in nonstandard models of arithmetic. This is
		conceptually distinct from the nonstandard methods used in
		\cite{HinmanKucaSchlesingerSheydvasser2019JNT}, where nonstandard extensions
		of $\N$ are employed as a tool inside ZFC to study the eventual behaviour of
		standard Ulam sequences.
		\item One could also ask for categorical or topos-theoretic formulations of
		Ulam rigidity, for instance by interpreting the Ulam construction in
		suitable arithmetic or realizability toposes and analysing the extent
		to which the rigidity phenomena are preserved.
		\item From the model-theoretic point of view, one might hope for analogues of the Marker–Steinhorn-type dichotomies \cite{MarkerSteinhorn1994DefinableTypes} for expansions of Presburger arithmetic by sparse additive sets such as Ulam sequences.
	\end{itemize}
	
	We leave these questions for future work. Even at the purely arithmetical
	level, there remain many open combinatorial problems about Ulam sequences; the
	results of this paper show that, whatever their ultimate resolution, the
	logical complexity of the corresponding statements is tightly controlled.
	
	% ------------------------------
	% Bibliography
	% ------------------------------

\end{document}